\documentclass[12pt]{amsart}
\usepackage{a4wide}
\usepackage{subfigure}
\usepackage{ifthen}
\usepackage{graphicx}

\newcommand\proofbox{\ensuremath{\blacksquare}\relax}
\newcommand{\noproofbox}{\ensuremath{\square}\relax}

\usepackage{amsmath}
\usepackage{amssymb}

\makeatletter
\def\swappedhead#1#2#3{%
  \thmnumber{\@upn{\the\thm@headfont #2\@ifnotempty{#1}{.~}}}%
  \thmname{#1}%
  \thmnote{ {\the\thm@notefont(#3)}}}
\makeatother

 \newtheoremstyle{changebreak}
   {9pt}
   {9pt}
   {\itshape}
   {}
   {\bfseries}
   {.}
   {\newline}
   {}

\renewenvironment{proof}[1][Proof]{\par\smallskip\noindent{\sc #1. }}{\qed}
\def\proofof #1 {\par\medskip\noindent {\sc Proof of #1. }}

\def\sketchof #1 {\par\medskip\noindent {\sc Sketch of proof of #1. }}

\def\qed{\rule{0pt}{0pt}\nolinebreak\hfill \proofbox \par\medskip}
\def\qedd{\rule{0pt}{0pt}\nolinebreak\hfill \noproofbox}

\newenvironment{remark}[1][]{%
  \par\noindent \textsc{\ifthenelse{\equal{#1}{}}{Remark. }{#1. }}}{%
   \par\medskip}

 \numberwithin{equation}{section}

\theoremstyle{changebreak}
\swapnumbers

\newtheorem{thm}{Theorem}[section]
\newtheorem{defn}[thm]{Definition}
\newtheorem{lem}[thm]{Lemma}
\newtheorem{cor}[thm]{Corollary}

\newtheorem{prop}[thm]{Proposition}
\newtheorem{conj}[thm]{Conjecture}
\newtheorem{question}[thm]{Question}
\newtheorem{observation}[thm]{Observation}

\newtheorem{deflem}[thm]{Definition and Lemma}

\swapnumbers

\newcommand{\ep}{\qed}

\newcommand{\PRconj}{\mathcal{G}}
\newcommand{\DRconj}{\mathfrak{g}}

\newcommand{\PRS}{G}
\newcommand{\PR}{G_{\s}}

\newcommand{\picturedir}{.}

\usepackage{hyperref}

\newcommand{\N}{\mathbb{N}}

\newcommand{\Z}{\mathbb{Z}}

\newcommand{\R}{\mathbb{R}}

\newcommand{\C}{\ensuremath{\mathbb{C}}}
\newcommand{\Ch}{\hat{\mathbb{C}}}

\newcommand{\D}{\mathbb{D}}
\newcommand{\Ds}{\D^*}

\newcommand{\eps}{\ensuremath{\varepsilon}}

\newcommand{\id}{\operatorname{id}}

\newcommand{\re}{\operatorname{Re}}
\newcommand{\im}{\operatorname{Im}}
\newcommand{\Log}{\operatorname{Log}}
\newcommand{\cl}[1]{\overline{#1}}

\newcommand{\Ek}{E_{\kappa}}

\newcommand{\ul}[1]{\underline{#1}}

\newcommand{\Sequ}{\mathcal{S}}

\newcommand{\Hplane}{\mathcal{H}}

\renewcommand{\u}{{\tt u}}

\newcommand{\bdyit}[2]
             {{\rule{0pt}{0pt}_{\mbox{$\scriptstyle #2$}}^{\mbox{%
                   $\scriptstyle #1$}} }}

\newcommand{\adds}{\ul{s}}
\newcommand{\s}{\adds}

\renewcommand{\r}{\ul{r}}

\newcommand{\extaddr}{\operatorname{addr}}
\newcommand{\itin}{\operatorname{itin}}

\newcommand{\gs}{g_{\adds}}

\newcommand{\periodic}[1]{\overline{#1}}
\newcommand{\per}[1]{\periodic{#1}}

\newcommand{\kappat}{\widetilde{\kappa}}
\newcommand{\wt}[1]{\widetilde{#1}}

\newcommand{\Et}{\widetilde{E}}

\newcommand{\F}{\mathcal{F}}
\newcommand{\ts}{t_{\s}}
\newcommand{\Xb}{\overline{X}}
\newcommand{\gammat}{\wt{\gamma}}

\newcommand{\compin}{\Subset}

\renewcommand{\P}{\mathcal{P}}

\newcommand{\J}{\mathcal{J}}
\newcommand{\I}{\mathcal{I}}
\newcommand{\B}{\mathcal{B}}

\begin{document}

\title[Topological Dynamics of Exponential Maps]{%
    Topological dynamics of exponential maps \\
       on their escaping sets}

\author{LASSE REMPE}
\address{Mathematics Institute \\University of Warwick 
  \\Coventry CV4 7AL\\ UK}
\email{lasse@maths.warwick.ac.uk}
\thanks{Supported in part
 by a postdoctoral fellowship of the 
 German Academic Exchange Service (DAAD).}

\begin{abstract}
 For the family of exponential maps $\Ek(z)=\exp(z)+\kappa$, we
  prove an analog of B\"ottcher's theorem by showing that 
  any two exponential maps $E_{\kappa_1}$ and $E_{\kappa_2}$ are
  conjugate on suitable subsets of their escaping sets, 
  and this conjugacy is quasiconformal. 
  Furthermore, we prove that any two attracting and
  parabolic exponential maps are conjugate on their sets of escaping
  points; in fact, we construct an 
  analog of Douady's
  ``pinched disk model'' for 
  the Julia sets of these maps. On the other hand, we show that two
  exponential maps are generally \emph{not} conjugate on their sets of
  escaping points.

 We also answer several
  questions about escaping endpoints of dynamic rays. In particular,
  we give a necessary and 
  sufficient
  condition for the ray to be continuously
  differentiable in such a point, and show that escaping points can
  escape arbitrarily slowly. Furthermore, we show that the principle of
  topological renormalization is false for attracting exponential
  maps. 
\end{abstract}

\maketitle


\section{Introduction}

 If $p$ is a polynomial of degree $d\geq 2$, then by B\"ottcher's Theorem 
  \cite[Theorem 9.1]{jackdynamics},
  $p$ is conjugate to 
  $z\mapsto z^d$ in a neighborhood of $\infty$. In the case where none of the
  critical points of $p$ is attracted to $\infty$ (or equivalently if
  the Julia set of $p$ is connected), this conjugacy can be extended
  to a biholomorphic mapping between the complement of the unit disk
  and the basin of infinity of $p$. The images of radial rays under
  this map give rise to the foliation of this basin by
  \emph{dynamic rays}, which have been used very successfully in the 
  combinatorial
  study of polynomials \cite{orsay}.

 In the family of exponential maps $\Ek: z\mapsto \exp(z)+\kappa$, the
  point $\infty$ is no longer an attracting fixed point, but rather an
  essential singularity, and the set of escaping points
   \[ I(\Ek):= 
     \{z\in\C:\Ek^n(z)\to\infty\} \]
  has no interior \cite[Section 2]{alexmisha} and is thus contained 
  in the Julia
  set. 
  Nevertheless, it was recently shown 
  by Schleicher and Zimmer
  \cite{expescaping} that $I(\Ek)$ is
  a union of curves to $\infty$
  which can be seen as an analog of
  dynamic rays of polynomials. 
  However, this still leaves
  open many questions
  on the topology of $I(\Ek)$, and on the dynamics
  of $\Ek$ thereon. For example, one can ask whether, as in the
  polynomial case, any two exponential maps with nonescaping singular
  values are conjugate on their sets of escaping points.

 We show, by a simple argument, that this is not true in general (see
  Section \ref{sec:noconjugacyexample}). In
  fact, it is already false when one of the parameters has an
  attracting fixed point and the other is a postsingularly finite (or
  \emph{Misiurewicz}) parameter; i.e., one for which the singular value
  $\kappa$ is preperiodic. (For polynomials,
  Misiurewicz parameters
  are among the most easily understood.)
  The argument generalizes to a much larger class of
  parameters; see Section \ref{sec:rigidity}.

 Despite these negative results, it is possible to make some statements 
  about the 
  topological dynamics on the set of escaping points in general. 
  We show the following, which can be seen as an analog of 
  B\"ottcher's theorem. (An \emph{escaping parameter} is one for which 
  the singular value escapes.) 

 \begin{thm}[Conjugacy Between Exponential Maps]    \label{thm:boettcheranalog}
  Let $\kappa_1,\kappa_2\in\C$. Let $R>0$ be large enough and consider
   the set
  \[ A := \{z\in\C: |E_{\kappa_1}^n(z)| \geq R \text{ for all $n\geq
             1$}\}. \]
  Then there exists a quasiconformal
   map $\phi:\C\to\C$ such that
  \begin{align}
    &\phi(E_{\kappa_1}(z)) =
       E_{\kappa_2}(\phi(z)) \text{ for all $z\in A$} \quad
             \text{and}  \label{eqn:conjugacy}
         \\
    &|E_{\kappa_1}^n(z) - E_{\kappa_2}^n(\Phi(z))| \to 0
                       \text{ for all $z\in I(E_{\kappa_1})\cap A$.}
     \label{eqn:samespeed}
  \end{align}

  If neither $\kappa_1$ nor $\kappa_2$ 
    are escaping
  parameters, then  
   $\phi|_{A\cap I(E_{\kappa_1})}$ 
   extends to a bijection 
     $\Phi: I(E_{\kappa_1})\to 
                     I(E_{\kappa_2})$ satisfying
   (\ref{eqn:conjugacy}) and (\ref{eqn:samespeed}) for all 
   $z\in I(E_{\kappa_1})$. 
 \end{thm}
 \begin{remark}[Remark 1]
  In general,
   the extended map $\Phi$ will not be 
   continuous.
 \end{remark}
 \begin{remark}[Remark 2]
  The number $R$ can be chosen of size $R=O(\max(|\kappa_1|,|\kappa_2|))$. 
 \end{remark}

 The proof of Theorem 
  \ref{thm:boettcheranalog} is achieved 
   by
  constructing an explicit model for the 
  set of escaping points, and
  then constructing a conjugacy 
  between an exponential
  map and this model on a suitable 
   set $A$ as in the
  theorem. In particular, this yields a 
  simplified 
  proof of the
  classification of escaping points given 
  by Schleicher and Zimmer
  (see Corollary \ref{cor:classification}). 
 
 While exponential 
  maps are generally not
  conjugate on their sets of escaping 
  points, the situation is quite
  different for parameters with an 
  attracting (or parabolic) periodic
  orbit, since such maps are 
  expanding on their Julia sets. 

 \begin{thm}[Topological Conjugacy] \label{thm:attractingconjugacy}
  Suppose that $\kappa_1$ and $\kappa_2$ are attracting (or parabolic)
   parameters. Then the map
    $\Phi$
  from Theorem \ref{thm:boettcheranalog} is a conjugacy.
 \end{thm}
 
 In fact, we give an explicit topological 
  model for the Julia set of
  such a parameter and the topological 
  dynamics thereon, 
  based on its combinatorics, as a 
  quotient of our
  general ``straight brush'' model. 
  Such models have already been
  constructed for the case of an attracting fixed point in
  \cite{aartsoversteegen} and for general periods in
  \cite{accessible}. (These constructions
  depended on the specific parameter, 
  and thus do not
  imply Theorem
  \ref{thm:attractingconjugacy}).

 On the other hand, we show a 
  somewhat surprising result for
  attracting exponential maps. For polynomials,
  \emph{renormalization} 
  fuels much of the detailed study of parameter
  spaces, ever since introduced by Douady and Hubbard. Although the
  concept does not generalize directly (due to the absence of
  compactness, there is no notion of ``exponential-like maps''), it
  was hoped 
  that some form of renormalization exists also in
  the exponential family (see \cite[Section VI.6]{habil}
  for a formulation of this question). 
  In particular, it was thought that
  renormalization is \emph{topologically} valid; i.e.~that, by
  collapsing certain rays and the regions between them for an
  exponential map $\Ek$ which is ``renormalizable'' of period $n$, the
  projection of the map $\Ek^n$ on this space will be conjugate to
  another exponential map. We show that this is false even for
  attracting exponential maps. More precisely, suppose that $\Ek$ has
  an attracting periodic orbit of period $n>1$, and let $E_{\kappa'}$
  be an exponential map with an attracting fixed point of the same
  multiplier as the attracting cycle of $\Ek$. If $U$ is an
  immediate attracting basin of $\Ek$, 
  then it is known (compare the
  discussion in Section \ref{sec:norenormalization}) that
   $\Ek$ restricted to $U$ is conformally
   conjugate to $E_{\kappa'}$ restricted to its Fatou set
   $F(E_{\kappa'})$. 

  \begin{thm}[No Topological Renormalization] \label{thm:renorm}
   The conformal conjugacy
        $\Psi:U \to F(E_{\kappa'})$
    does not extend 
    continuously to $\partial U$.
  \end{thm}

 The classification of escaping points by Schleicher and Zimmer
  \cite{expescaping} exposed a feature of dynamic rays which does not
  occur for polynomials: some dynamic rays have endpoints which also
  escape to $\infty$. For the purposes of their
  classification, Schleicher
  and Zimmer used topological arguments to reach these endpoints and
  do not provide much 
  additional information about them. In particular, it is not
  a priori
  clear whether the escape speed of these endpoints is independent of
  the parameter. That this is the case follows from
  Theorem \ref{thm:boettcheranalog}. Because our model for escaping
  points is very explicit, we can also answer several other questions
  concerning these endpoints; in particular we show that they can
  escape arbitrarily slowly, or in fact with any prescribed
  escape
  speed. (Escaping points which are not endpoints
  are known to always escape with iterated exponential speed.)
 \begin{thm}[Arbitrary Escape Speed] \label{thm:arbitrarilyslowescape}
  Let $\kappa\in\C$. 
  Suppose that $r_n$ is a sequence of positive real numbers such that
  $r_n\to \infty$ and $r_{n+1}\leq \exp(r_n) + c$ for some $c>0$. Then
  there is an escaping endpoint
  $z\in I(\Ek)$ and some $n_0\in\N$ such that, for
  $n\geq n_0$,
   $|\re(\Ek^{n-1} (z)) - r_n | \leq 2+2\pi$.
 \end{thm}
 
Furthermore, we give an explicit
  necessary and sufficient condition (independent of the parameter)
  under which a ray with escaping endpoint is differentiable
  in this endpoint (Theorem \ref{thm:endpointdifferentiability}). 

 Finally, we also discuss the situation in parameter space. The
  set $\I$ of parameters for which the singular value lies on a dynamic
  ray can be described in terms of 
  \emph{parameter rays} \cite{markusdierk}; 
  this result is extended to
  escaping endpoints in \cite{markuslassedierk}.
  As in the
  dynamical plane, it is interesting to ask what the topology of this
  set is. In Section \ref{sec:parameterspace}, we show that the
  bijection between our model space and the set of escaping parameters
  is a homeomorphism on large sets, and that there exists a sequence
  of ``Cantor Bouquets'' in parameter space whose union covers $\I$.

\subsection*{Organization of the article}
 After a short exposition of an
  example of two exponential maps which are not conjugate on their
  escaping sets in Section \ref{sec:noconjugacyexample}, the
  fundamental construction of our model and the correspondence with
  the set of escaping points of an exponential map 
  (including the proof of Theorem
   \ref{thm:boettcheranalog}) are carried out in
  Sections \ref{sec:model} and \ref{sec:classification}. These two
  sections form the backbone of the remainder of the article.

 The remaining sections are mostly independent of each other.
  Section \ref{sec:limitsofrays} contains a short discussion
  of the limiting
  behavior of dynamic rays, collecting results and techniques
  from \cite{expescaping} and \cite{markuslassedierk}.
  In Section \ref{sec:differentiability} we
  answer the question which rays are differentiable in their escaping
  endpoints and in Section \ref{sec:speedofescape}, we prove some
  facts on the speed of escape (including
  Theorem \ref{thm:arbitrarilyslowescape}), demonstrating that questions
  of this type can easily be answered using our model.

 Section \ref{sec:rigidity} discusses the uniqueness of our
  correspondence and some further results on escaping set rigidity of
  exponential maps. Sections \ref{sec:topologyattracting} 
  and \ref{sec:norenormalization}
  are devoted to the proofs of Theorems \ref{thm:attractingconjugacy} and
   \ref{thm:renorm}, respectively.
  Finally, Section \ref{sec:parameterspace}
  examines continuity properties of parameter rays and Section
  \ref{sec:openquestions} discusses some questions that remain open.

\subsection*{Remarks on notation}
  The Julia, Fatou and escaping sets of an exponential map will be
   denoted $J(\Ek)$, $F(\Ek)$ and $I(\Ek)$, as usual. It is 
   well-known \cite{bakerexp,alexmisha} that exponential maps
   have no wandering domains and at most one nonrepelling cycle. 

  In particular, if $F(\Ek)\neq \emptyset$, then
   $F(\Ek)$ consists of the basin of attraction of some attracting or
   parabolic cycle, or of the iterated preimages of a Siegel disk. 
   In this case, 
    the map $\Ek$ (and also the parameter $\kappa$) is called 
    \emph{attracting}, \emph{parabolic} or \emph{Siegel}, respectively. 
   If the \emph{postsingular set} 
     \[ \P(\Ek) := \overline{\bigcup_{n\geq 0} \Ek^n(\kappa)} \]
   is finite, then $\kappa$ is called a \emph{Misiurewicz} parameter;
   in this case, $J(\Ek)=\C$. 

  Throughout the article, we shall fix the function
   $F:[0,\infty)\to[0,\infty);t\mapsto \exp(t)-1$
   as a model for exponential growth. We shall routinely make use of the
   fact that, for all $t_2 \geq t_1 \geq 0$, 
  \begin{equation}
   F(t_2) - F(t_1) =
      e^{t_1} F(t_2 - t_1) \geq
      F(t_2 - t_1)  \label{eqn:propertiesofF}
  \end{equation}

  Following \cite{expescaping},
   a sequence $\s = s_1 s_2 s_3 \dots \in \Z^{\N}$ 
   of integers is called 
   an \emph{external address}; the shift map on external addresses
   is denoted by $\sigma$. We end any proof by the symbol
   $\proofbox$; statements which are cited without proof
   are concluded by $\noproofbox$.

  \section{Two Exponential Maps not Conjugate on Their Escaping Sets}

 \label{sec:noconjugacyexample}

 In this section, we describe, as a motivation for our further
 discussions, an example of two exponential maps which are not
 conjugate on their sets of escaping points. The two maps we will discuss
 are fairly well-understood since the 1980s.

 Set $\kappa := -2$ and $E:=E_{\kappa} : z\mapsto \exp(z)-2$.
  This
  map has an attracting fixed point between $-2$ and $-1$ which
  attracts the entire interval $(-\infty,0)$. In
  particular, the singular value $\kappa$ is contained in the Fatou
  set. It is well known \cite[Proposition 3.3]{devgoldberg} that the
  Julia set of $E$ is a ``Cantor Bouquet''; in particular every connected
  component is an injective curve $\gamma:[0,\infty)\to\C$, where
  $\gamma(t)$ escapes for $t>0$;
  i.e. $\gamma\bigl((0,\infty)\bigr)
                               \subset I(E)$.

 On the other hand, let $\kappat:= \log(2\pi)+\frac{\pi}{2}i$, and denote
   \[ \Et := E_{\kappat}: z\mapsto \exp(z) +\log(2\pi)+\frac{\pi}{2}i. \]
  Note that $\Et(\kappat) = 2\pi i + \kappat$ is a fixed point. The
   Julia set of $\Et$ is the entire plane, and the dynamics of $\Et$
   is less well-understood
   than than that of $E$. Nonetheless, it has long been known
   that there exists an injective curve of escaping points
   landing at the fixed point $\Et(\kappat)$ and tending to $\infty$
   in the other direction. (For a proof, see
   \cite[Theorem 3.9]{dghnew1} or \cite[Proposition 6.11]{expescaping}; for
   arbitrary Misiurewicz parameters the same fact is proved in
   \cite[Theorem 4.3]{expper}.)
   Pulling back, 
   we obtain a curve
   $\gammat : [0,\infty)\to \C$
   with $\gammat(0)=\kappat$, $\gammat\bigl((0,\infty)\bigr)\subset
   I(E_{\kappat})$ and $\gammat(t)\to\infty$ for $t\to\infty$.

In the following, we will denote the Julia set of $E$ by $J$ and that
of $\Et$ by $\wt{J}$; similarly for their sets of escaping points etc.
We are now ready to prove that these two maps are not conjugate on
their sets of escaping points.

\begin{prop}[No Conjugacy] \label{prop:noconjugacy}
 The maps $E_{|I}$ and $\Et_{|\wt{I}}$ are not topologically conjugate.
\end{prop}
\begin{proof}
 Assume, by contradiction,
 that there is a conjugacy between 
  $E$ and $\Et$ on their
 sets of escaping points, say 
  $\Phi:I\to\wt{I}$.
 Consider the curve $\gamma$ 
  corresponding to $\gammat$ under 
  $\Phi$, i.e.~%
 $\gamma:=\Phi^{-1}\circ 
        \gammat:(0,\infty)\to I$. 
  We first claim
 that 
 $\lim_{t\to\infty} \gamma(t)=\infty$. 
  Suppose not; by the above description
 of the topology of $J$ we would then 
  have $\lim_{t\to\infty}\gamma(t)=z$ 
  for
 some $z\in\C$. However, then $E(z)$ 
  would be a fixed point of $E$, and
  all points
  on $\gamma$ converge to $E(z)$ 
  under iteration. This is a
 contradiction, as the map $E$ has only 
  repelling periodic points in
 its Julia set.

 Again by the topology of $J$, 
  $\gamma$ has an endpoint 
 $z_0:=\gamma(0) := 
   \lim_{t\to 0} \gamma(t)$. Pick any 
  point $w\in I$, say
  $w=2$, and denote 
  $\wt{w} := \Phi(w)$. Choose any 
  open neighborhood $U$
 of $w$ such that 
    \[ \Phi(U\cap I) \subset 
        \{z\in\C:|z-\wt{w}|<1\}. \]
 Now, since $w\in J$,
 we can find a point $z_1\in U$ with 
  $E^n(z_1)=z_0$ for some $n$.
 Pulling back $\gamma$ along the 
  corresponding branch of $E^{-n}$, 
  we obtain
 a curve $\alpha:(0,\infty)\to I$ with 
  $\lim_{t\to 0} \alpha(t) = z_1$. In
 particular, there exists $t_0>0$ with 
   $\alpha(t)\in U$ for $t\leq t_0$. 

 Now consider the curve 
   $\wt{\alpha}:=\Phi\circ \alpha$. 
   This curve satisfies
   \[ |\wt{\alpha}(t)-\wt{w}|<1 \]
  for $t\leq t_0$. On the other hand,
   \[ \lim_{t\to 0}\Et^n(\wt{\alpha}(t)) 
     = 
    \lim_{t\to 0}\gammat(t) = \kappat, 
        \]
  and thus 
  $\lim_{t\to 0} |\wt{\alpha}(t)|=
                         \infty$. 
  This is a
   contradiction. 
\end{proof}

For further discussion and more general 
 results 
 on the (non-)existence of conjugacies 
 on the sets of
 escaping points, see Section 
 \ref{sec:rigidity}.

\section{A Model for the Set of Escaping Points} \label{sec:model}

To motivate the definitions which follow,
 let us shortly review the structure of the Julia set for the map $E=E_{-2}$
 of the previous section,
 as described e.g.~in \cite{devgoldberg}. Because the interval
 $(-\infty,-2)$ is contained in the Fatou set, the Julia set $J$ is
 disjoint from the preimages of this interval, which are straight
 lines of the form $\{\im z = (2k-1)\pi\}$. In other words, $J$ is completely
 contained in the strips 
   \[ S_k := \Bigl\{z:\im z \in
                    \bigl((2k-1)\pi,(2k+1)\pi\bigr)\Bigr\}. \]
 To any point $z\in J$ we can thus associate an external address $\s$
 such that $\Ek^{k-1}(z)\in S_{s_k}$ for all
 $k$. This sequence is called the \emph{external address of $z$}.
 It turns out that the connected components of $J$ are exactly the
 sets of the form
  \[ \{z\in\C: z \text{ has external address $\s$} \}, \]
 which are curves consisting of escaping points together with an
 endpoint which may or may not escape.

Let us now develop the promised model for the set of escaping points of an
  exponential map (with nonescaping singular orbit).
  Based on the description above, our model should consist of
  pairs $(\s,t)\in \Z^{\N}\times [0,\infty)$\footnote{%
   Note that not all external addresses can be realized by an
    exponential map; accordingly the same will be true of our model.}).
  Note that the space of external addresses has a natural
  topological structure, namely that induced by
  the lexicographic order on external
  addresses (open sets are unions of open
  intervals). Thus $\Z^{\N}\times [0,\infty)$ 
  is equipped with the product of
  this topology and the
  usual topology of the real numbers. 

 For a given point $(\s,t)$, the first entry $s_1$ of $\s$ should be thought
  of as the imaginary part, while $t$ corresponds to the
  real part.
  We thus define $Z(\s,t):=t+2\pi i s_1$
  and abbreviate
  $|(\s,t)|:=|Z(\s,t)|$. We shall also write
  $T$ for the projection to the second component; i.e. $T(\s,t)=t$. In
  analogy to the potential-theoretic interpretation of dynamic rays
  of polynomials, we will sometimes refer to $t$ as the ``potential''
  of the point $(\s,t)$.

 We now define a model function which will naturally give rise
  to our model space. There is
  considerable freedom in the definition; to suit our needs,
  we have chosen here to use a 
  function which allows very explicit calculations. 

 Our model dynamics is 
  then given by
     \[ \F(\s,t) := (\sigma(\s),F(t)-2\pi |s_2|)). \]
 Its key feature is that, as for exponential maps,
  the size of the image of a point is
  roughly the exponential of its real part. Indeed, 
  if  $\s\in \Z^{\N}$
  and $t\geq 0$ with $T(\F(\s,t))\geq 0$, then 
  \[ \frac{1}{\sqrt{2}}F(t) \leq |\F(\s,t)| \leq F(t).\]
  We now define
  \begin{align*}
   \Xb &:= \{(\s,t): \forall n\geq 0, T(\F^n(\s,t))\geq 0\} \ \text{
  and} \\ 
    X &:= \{(\s,t)\in \Xb: T(\F^n(\s,t))\to\infty\}.
 \end{align*}
 
 The space $X$ will be our model of the set of escaping
 points; as we show in Section
 \ref{sec:topologyattracting}, $\F|_{\Xb}$ is conjugate to
 the attracting exponential map $E$
 considered in Section \ref{sec:noconjugacyexample} on its Julia set. In
 particular, the set $\Xb$ is a ``straight brush'' in the sense of
 \cite{aartsoversteegen}, which we show directly in the following
 observation. 

 \begin{observation}[Comb Structure of $\Xb$ and $X$]
  For every external address $\s$, there exists a $t_{\s}$, 
  $0\leq \ts \leq \infty$, such that
   \begin{equation*}
      \{t\geq 0: (\s,t)\in \Xb\} = [\ts , \infty); 
   \end{equation*}
  this $t_{\s}$ depends lower semicontinuously on $\s$. Furthermore, 
   \begin{align}
     (\ts,\infty) &\subset X_{\s} := \{t\geq 0: (\s,t)\in X\} \quad \text{and}
         \label{eqn:Xs} \\
      \F(\s,\ts) &= (\sigma(\s), t_{\sigma(\s)}). 
        \label{eqn:tsfunctionalequation}
   \end{align}
 \end{observation}
 \begin{proof}
  Suppose that $(\s,t)\in \Xb$ and $t'=t+\delta$, $\delta>0$. By the
  definition of $\F$, we have
  \[
     T(\F(\s,t')) - T(\F(\s,t)) =
     F(t')-F(t) \geq F(t' - t) = F(\delta). \]
  By induction,
  \[ T(\F^n(\s,t')) \geq T(\F^n(\s,t')) - T(\F^n(\s,t)) \geq
       F^n(\delta)\to \infty. \]
  This proves the first claim as well as
   (\ref{eqn:Xs}). Note also that
   (\ref{eqn:tsfunctionalequation}) follows directly
   from the definitions. To prove semicontinuity, note that
  $\Xb$ is a closed set. Therefore, for any $R>0$ the set
   $ \{\s:t_{\s} \leq R \} =
      \{\s:(\s,R)\in \Xb\}$
  is closed. 
 \end{proof}

Following \cite{expescaping}, we will call an external address $\s$
  \emph{exponentially bounded} if $\ts <\infty$. Furthermore, we will
  call such an address \emph{fast} if
  $(\s,\ts)\in X$; otherwise $\s$ is called \emph{slow}. (It is not
  difficult
  to see that these definitions are equivalent to those given
  in \cite{expescaping}, compare Corollary \ref{cor:properties}.)
  We will also denote the space of all exponentially bounded
  external addresses by $\Sequ_0$.


\section{Classification of Escaping Points}
\label{sec:classification}

\begin{figure}
 \includegraphics[width=\textwidth]{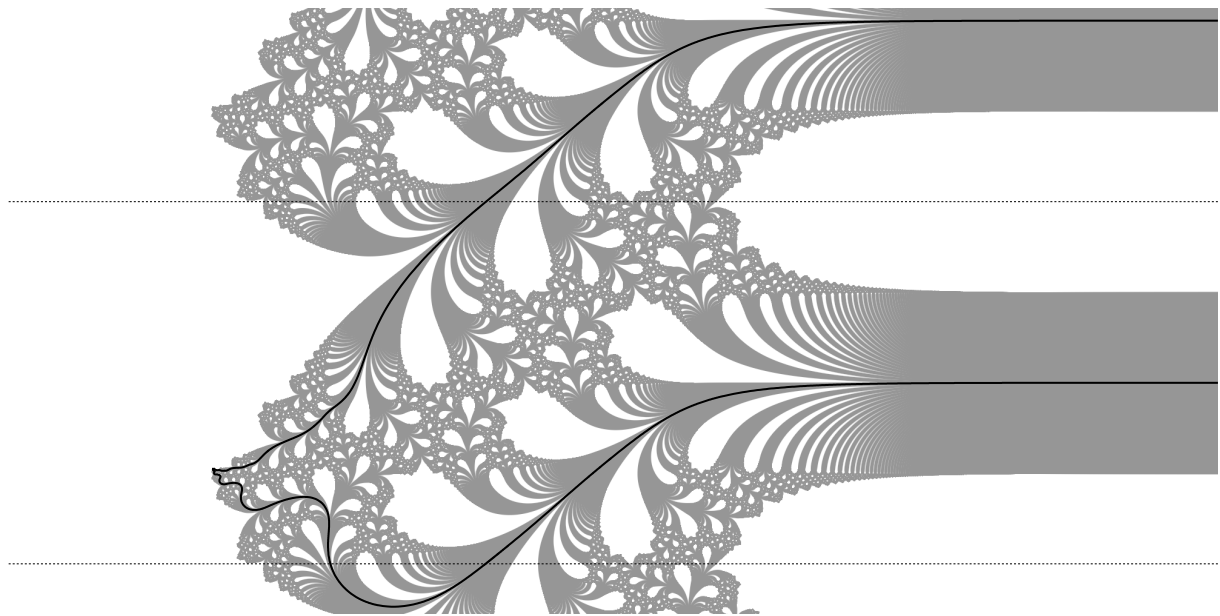}
 \caption{Dynamic rays may cross strip boundaries.\label{fig:rayscrossingboundaries}}
\end{figure}

 In the following, we fix 
   some arbitrary exponential map
   $\kappa\in\C$. As before, we say 
   that a point $z\in\C$ has external
   address $\s$ if
   \[ \im(\Ek^{n-1}(z))\in S_{s_n} \]
   for all $n$, where 
     $S_k = \bigl\{z:\im z \in \bigl((2k-1)\pi,(2k+1)\pi\bigr)\bigr\}$.
  Note that, for general $\kappa$, not all points
  $z\in I(\Ek)$ have an external address, as components of $I(\Ek)$
  may cross the strip boundaries (see Figure
  \ref{fig:rayscrossingboundaries}). However, some forward iterate of
  $z$ always has an external address. Indeed, 
   $|\Ek^n(z) - \kappa | = \exp( \re \Ek^{n-1}(z))$,
  so $\re \Ek^{n-1}(z)\to\infty$. Thus, if $n$ is large enough, the
  orbit of $\Ek^n(z)$ will be contained in the half plane 
  $\{\re z > \re\kappa\}$. Images of points on the strip boundaries,
  on the other hand, lie in $\{\kappa - t : t > 0\}$, so the orbit of
  $\Ek^{n-1}(z)$ never intersects these boundaries and thus has an
  external address.

  We now construct a
  conjugacy between $\F$ and $\Ek$ (defined on a suitable subset of $X$). 
  This is done by iterating forward in our model and
  then backwards in the dynamics of $\Ek$. To this end, we define the
  inverse branches $L_k$ of $\Ek$ by
  \[ L_k(w):=\Log(w-\kappa) + 2\pi i k, \]
  where 
  $\Log:\C\setminus (-\infty,0]\to S_0$ 
  is the principal branch
  of the logarithm. Thus $L_k(w)$ 
  is defined and analytic whenever $w-\kappa\notin
  (-\infty,0]$. 

 Define maps $\DRconj_k$ inductively by $\DRconj_0(\s,t):=Z(\s,t)$ and
  \[ \DRconj_{k+1}(\s,t):= L_{s_1}(\DRconj_k(\F(\s,t))) \]
 (wherever this is defined).

 Fix $K>2\pi+6$ such that $|\kappa|\leq K$, and let us define
 \begin{align}
    \label{eqn:Q}
  Q &:= Q(K) := \max\{\log(4(K+\pi+3)),\pi+2\} \quad\text{and} \\
  \notag
  Y &:= Y_Q := 
        \left\{(\s,t)\in\Xb:T(\F^n(\s,t))\geq Q \text{ for all $n$} \right\}.
 \end{align}

 Note that $Y$ contains the set
  $\{(\s,t): \s\in\Sequ_0\text{ and }t\geq \ts+Q\}$. Note also that, 
  for every
  $x=(\s,t)\in X$, there exists some $n$ such that $\F^n(x)\in Y$.

 \begin{lem}[Bound on $\DRconj_k$] \label{lem:boundeddistance}
  For all $k$, the map $\DRconj_k$ is defined on $Y$.
   For all $(\s,t)\in Y$, we have
   $|\re \DRconj_k(\s,t) - t|<2$, and in particular
   $|\DRconj_k(\s,t)-Z(\s,t)|<\pi + 2$.
 \end{lem}
 \begin{proof}
  The idea of the proof, which proceeds by induction, is
  quite simple. By the induction hypothesis, we know that
  $\DRconj_k(\F(\s,t))$ and $Z(\F(\s,t))$ are close, and by the definition
  of $\F$, the values $|Z(\F(\s,t))|$ and $|\Ek(Z(\s,t))|$ are
  essentially the same, namely $F(t)$ (up to a constant
  factor). Pulling back $\DRconj_k(\F(\s,t))$ and $\Ek(Z(\s,t))$ by the same
  branch of $\Ek^{-1} = \Log(z-\kappa)$, this constant factor
  translates to an additive constant for the real parts, as
  desired. The somewhat unpleasant calculations 
  which follow flesh out
  this idea and fix the constants.

  To begin the induction, note that
  the case $k=0$ is trivial. Now fix $k\geq 0$ such that the
  claim is true for $k$; we will show
  that it holds also for $k+1$. Let $(\s,t)\in Y$. By the induction
  hypothesis,
  \begin{align*}
     |\DRconj_k(\F(\s,t))-Z(\F(\s,t))| &\leq \pi + 2 \quad \text{and} \\
     \re\left(\DRconj_k(\F(\s,t))\right) 
      \geq T(\F(\s,t)) - 2 \geq
            Q - 2  &\geq 0. 
  \end{align*}
  Furthermore, we have
  \[ F(t)=\exp(t)-1 \geq \exp(Q) -1 \geq 4(K+\pi + 2) \]
  In particular,
  \[ |\DRconj_k(\F(\s,t))| \geq
     |Z(\F(\s,t))| - \pi - 2
     \geq \frac{F(t)}{\sqrt{2}} - \pi - 2 > 2K. \]
  Thus $\DRconj_k(\F(\s,t)) - \kappa \notin (-\infty, 0]$, so $\DRconj_{k+1}(\s,t)$
  is defined. Furthermore, we can write
  \[ \DRconj_k(\F(\s,t))-\kappa =
     Z(\F(\s,t)) + \bigl(\DRconj_k(\F(\s,t))-Z(\F(\s,t))-\kappa\bigr), \]
  and, by the definition of $Q$,
   \[ |\DRconj_k(\F(\s,t))-Z(\F(\s,t))-\kappa| \leq
       \pi + 2 + K \leq \frac{1}{4}\exp(Q)-1\leq \frac{1}{4}\exp(t) - 1. \]

  Therefore
   \begin{align*}
      \re(\DRconj_{k+1}(\s,t)) &= \log|\DRconj_k(\F(\s,t)) - \kappa | \geq 
      \log\left( |Z(\F(\s,t))| + 1 - \frac{1}{4}\exp(t)\right) \\
      &\geq
      \log\left( \frac{1}{\sqrt{2}} \exp(t)-\frac{1}{4}\exp(t) \right) >
      \log\left(\frac{\exp(t)}{4}\right) =
     t - \log 4. 
   \end{align*}
  Analogously $\re(\DRconj_{k+1}(\s,t))\leq t + \log 4$, and thus
   $|\re(\DRconj_{k+1}(\s,t) - t | \leq \log 4 < 2$, as required.
 \end{proof}

 With the estimate of Lemma \ref{lem:boundeddistance}, we can now
  construct the required conjugacy by a standard contraction argument.

  \begin{thm}[Convergence to a 
                          Conjugacy] 
      \label{thm:conjugacyconvergence}
   On $Y$, the functions $\DRconj_k$ 
   converge uniformly (in $(\s,t)$ and
   $\kappa$ with $|\kappa|\leq K$) to a function
   $\DRconj:Y\to J(\Ek)$ such that $\DRconj(\s,t)$ has
   external address $\s$ for each $(\s,t)\in Y$. This function
   satisfies ${\DRconj}\circ \F = \Ek \circ {\DRconj}$ and 
   \begin{equation}
      \Bigl| \DRconj(\s,t) - \bigl(\Log(Z(\F(\s,t)))+2\pi i s_1 \bigr) \Bigr| \leq
        e^{-t}\cdot (2K+2\pi +4). \label{eqn:asymptotics}
   \end{equation}

   Furthermore, $\DRconj(\s,t)\in I(\Ek)$ if and only if
    $(\s,t)\in X$; ${\DRconj}$ is
   a homeomorphism between $Y$ and its image; and 
   $\DRconj(\s,t)$ depends holomorphically on $\kappa$ for fixed $(\s,t)\in
    Y$. 
  \end{thm}
 \begin{remark} Note that, for every $\s\in\Sequ_0$,
    (\ref{eqn:asymptotics}) implies that
    \begin{equation} \DRconj(\s,t) = t + 2\pi i s_1 + O(e^{-t}). \label{eqn:rayasymptotics}
                                            \end{equation}
 \end{remark}
 \begin{proof}
   Recall from the previous proof that, for $n\geq 1$,
   \begin{align*} 
   &|\DRconj_k(\F^n(\s,t))|>2K \quad (\ \geq |\kappa| + 2\pi + 6\ ),
       \quad \text{and} \\
   &\re\left(\DRconj_k(\F^n(\s,t))\right)>0.
   \end{align*}
   Furthermore, the distance between $\DRconj_k(\F(\s,t))$ and
   $\DRconj_{k+1}(\F(\s,t))$ is at most $2\pi+4$. Thus we can connect these
    two points by a straight line within the set
   \[ \{z\in\C: |z-\kappa|\geq 2 \text{ and } z-\kappa\notin\ 
                (-\infty,0]\}.\]
   Since $L_{s_1}'(z)=
              \frac{1}{z-\kappa}$,
   \[ |\DRconj_{k+1}(\s,t)-\DRconj_{k+2}(\s,t)| \leq 
      \frac{1}{2} |\DRconj_{k}(\F(\s,t))-\DRconj_{k+1}(\F(\s,t))|. \]
   It follows by induction that
   \[ |\DRconj_{k+1}(\s,t)-\DRconj_{k+2}(\s,t)| \leq
       2^{-(k+1)}|Z(\F^{k}(\s,t))-
                        \DRconj_1(\F^k(\s,t))|\leq
       2^{-(k+1)}(\pi+2), \]
   so the $\DRconj_n$ converge 
   uniformly on $Y$. By definition,
   $\DRconj(\s,t) = L_{s_1}(\DRconj(\F(\s,t)))$,
   and thus
    \[ \Ek \circ G = G \circ \F. \]

\smallskip

   To prove the asymptotics 
    (\ref{eqn:asymptotics}), 
     we first observe that
 \[     |Z(\F(\s,t))| \geq
              F(t)/\sqrt{2} \geq
              2\exp(t)/3  \]
     (which follows easily
      from the fact that
     $t \geq 3$), and thus
  \[
   \left|\frac{\DRconj(\F(\s,t))-\kappa-
                        Z(\F(\s,t))}{%
                                Z(\F(\s,t))}\right|
    \leq \frac{3(\pi+2+K)}{%
                          2\exp(t)}
    \leq \frac{3}{8}.
  \]
 Since
       $|\Log(1+z)| \leq 4|z|/3$ when
        $|z|\leq 3/8$, it follows that
    \begin{align*}
       \Bigl| \DRconj(\s,t) - 
               \bigl(\Log(Z(\F(\s,t)))
                       +2\pi i s_1 \bigr) \Bigr|
     &= \left| \Log\left(1+
            \frac{\DRconj(\F(\s,t))-\kappa - 
                               Z(\F(\s,t))}{%
                        Z(\F(\s,t))}\right)\right| \\
     &\leq  
      \frac{4}{3} \frac{3(\pi+2+K)}{2\exp(t)}
     = 
      e^{-t}\cdot(2\pi+4+2K).    
    \end{align*}

 \smallskip

   By
   Lemma \ref{lem:boundeddistance}, $\DRconj(\s,t)$ escapes under
   iteration of $\Ek$ if
   and only if $(\s,t)$ escapes under iteration of $\F$.
   Clearly the point $\DRconj(\s,t)$ has the correct external address (note
   that $\arg \bigl(\DRconj_{k-1}(\s,t)\bigr)$ is bounded away from
   $\pm \pi$, so that the values $\DRconj_k(\s,t)$ cannot converge to the strip boundaries).
   In
   particular, $\DRconj(\s,t)\neq \DRconj(\s',t')$ whenever $\s\neq \s'$, because
   the points have different external 
   addresses. On the other hand,
   the $\F$-orbits
   of $(\s,t)$ and $(\s,t')$, $t\neq
   t'$, will eventually be arbitrarily far 
   apart, and therefore the
   same holds for $\DRconj(\s,t)$ and 
   $\DRconj(\s,t')$ under $\Ek$. This proves
   injectivity.

  The function ${\DRconj}$ is continuous as uniform limit of continuous
  functions; for the same reason, $\DRconj(\s,t)$ is analytic in
  $\kappa$. To prove that the inverse $\DRconj^{-1}$ is continuous, note
  that we can compactify both $Y$ and $\DRconj(Y)$ by adding a point at
  infinity. The extended map ${\DRconj}$ is still continuous,
  and the inverse of a continuous
  bijective map on a compact space is continuous. 
 \end{proof}

\begin{remark}
 The asymptotic description of $\DRconj(\s,t)$ in terms of
  \begin{align*}
     \Log( Z(\F(\s,t))) + 2\pi i s_1 = &
      \log\sqrt{(F(t)-2\pi |s_2|)^2 + (2\pi s_2)^2}\\
      &+i\arg\bigl((F(t)-2\pi|s_2| + 2\pi i s_2\bigr) + 2\pi i s_1 
  \end{align*}
 is somewhat awkward. Had we used, instead
of $\F$, the map
  \[ \F'(\s,t) := (\sigma(\s),\sqrt{F(t)^2 - (2\pi s_2)^2}), \]
then the whole construction would have carried through analogously 
 (with somewhat improved
  constants). For the map 
  ${\DRconj}':Y\to J(\Ek)$ that
  we obtain this way, we would correspondingly have the following asymptotics:
\[ |\re(G'(\s,t)) - t | < e^{-t}\cdot (2K+2\pi+4). \]
  In this article, we will never use the asymptotics in any
  other form than (\ref{eqn:rayasymptotics}), whereas
  we shall rather
  often make direct calculations. This is why we have opted to
  use the function $\F$ rather than $\F'$. 
\end{remark}

In order to extend ${\DRconj}$ to a bijection
 $\DRconj:X\to I(\Ek)$, 
  the main remaining problem is to decide when a point is
 contained in $\DRconj(Y)$. The
 following is a counterpart to Theorem
 \ref{thm:conjugacyconvergence}. 

 \begin{thm}[Points in the 
                     Image of ${\DRconj}$] 
    \label{thm:surjectivity}
  Suppose that $z\in\C$ spends its entire orbit in the halfplane
   $\{w\in\C:\re w \geq Q+1 \}$. Then $z$ has an external address 
   $\s$, and
   there exists $t\geq \ts$ 
   such that $(\s,t)\in Y$ and $z=\DRconj(\s,t)$.
 \end{thm}
 \begin{proof} 
  First note that, for $n\geq 1$, 
   \[ |\Ek^n(z)-\kappa| = 
     \exp(\re(\Ek^{n-1}(z))) \geq 
         \exp(Q+1) > 2K \]
  and $\re(\Ek^n(z))>0$, so 
  $\Ek^n(z)-\kappa\notin (-\infty,0)$. Therefore, no
  iterate of $z$ lies on the strip boundaries, and thus $z$
  has an external address $\s$. 

  Consider the sequence $t_k$ of potentials uniquely defined by
   \[ \F^k(\s,t_k) = (\sigma^k(\s),\re(\Ek^k(z))). \]
   We claim that
    (similarly to Lemma
     \ref{lem:boundeddistance}), for
     $j\leq k$,
   \begin{equation}
     \left| T(\F^j(\s,t_k)) - 
      \re(\Ek^j(z)) \right| \quad\left(\  =
     \left| \log\frac{T(\F^{j+1}(\s,t_k))
                          + 2\pi s_{j+2}+1}{%
                  |\Ek^{j+1}(z)-\kappa|}
       \right|\  \right) \quad
         \leq 1.   \label{eqn:pointsinimage}
   \end{equation}
   The idea is again to prove this
     by induction:
     since 
     $\Ek^{j+1}(z)$ belongs to
     the strip $S_{s_{j+2}}$ and
     $T(\F^{j+1})(\s,t_k)$ 
     and $\re\Ek^{j+1}$ are close
     by the induction hypothesis, the
     ratio in the logarithm is bounded.
     We omit the precise calculations
     here.

  Let $t$ be any limit point of the sequence
   $t_k$. Then
   $\left| T(\F^j(\s,t)) - \re(\Ek^j(z))
                           \right| \leq 1$
   by (\ref{eqn:pointsinimage});
   in particular, $(\s,t)\in Y$.

   Since $\DRconj(\s,t)$ also has external 
     address $\s$, it now follows that
     the distance between 
     $\Ek^j(\DRconj(\s,t))$ and $\Ek^j(z)$ 
     is bounded
   for all $j$. By the same contraction argument as in the proof of Theorem
   \ref{thm:conjugacyconvergence}, they are equal.
 \end{proof}

We can now prove the existence of a global correspondence between $X$
and $I(\Ek)$.

\begin{cor}[Global Correspondence] \label{cor:globalcorrespondence}
 Suppose that $\kappa\notin I(\Ek)$. Then
 ${\DRconj}|(Y\cap X)$ extends to a bijective function
   \[ G : X \to I(\Ek) \]
 which satisfies $\DRconj(\F(\s,t)) = \Ek(\DRconj(\s,t))$. 
 The function ${\DRconj}$ is a homeomorphism
   on every $\F^{-k}(Y\cap X)$ and for every $\s\in\Sequ_0$
   the function $t\mapsto \DRconj(\s,t)$ is continuous (``${\DRconj}$ is continuous
   along rays'').

 If $\kappa\in I(\Ek)$, then ${\DRconj}$ has a similar extension as follows. 
   There exists
 $(\s^0,t^0)\in X$ with $\DRconj(\s^0,t^0)=\kappa$, and ${\DRconj}$ is 
 defined for all $(\s,t)\in X$ except
  those with $\F^n(\s,t)=(\s^0,t')$ for some $n\geq 1$ and
  $t'\leq t^0$. For every $k$, ${\DRconj}$ is
  continuous on the intersection of $\F^{-k}(Y\cap X)$ with its domain
  of definition, and it is continuous along rays.

 Whenever $x_0\in \F^{-n}(Y)$ such that $\DRconj(x_0)$ is defined, there exists
  a neighborhood $U$ of $x_0$ in $\F^{-n}(Y)$ such that ${\DRconj}|_U$ is
  defined and continuous (as a function of $x=(\s,t)$ and $\kappa$) for
  all parameters with $|\kappa|\leq K$.
\end{cor}
\begin{proof} We will consider only the case of $\kappa\notin I(\Ek)$; the
 other statements follow similarly.

 It is sufficient to show, by induction, that ${\DRconj}$
  extends to a homeomorphism
   \[\DRconj:\F^{-k}(Y\cap X) \to \Ek^{-k}(\DRconj(Y\cap X))\]
  for every $k\geq 1$.
  Indeed, the
  sets of definition 
  clearly exhaust all of $X$, while the range exhausts $I(\Ek)$ by Theorem
   \ref{thm:surjectivity}. Continuity along rays also follows
  because every $(\s,t)\in X$ has a neighborhood on the ray that is
  completely contained in the same $\F^{-k}(Y\cap X)$.

 So let us suppose that ${\DRconj}$ has been 
   extended to $\F^{-k}(Y\cap X)$. 
  First note that we can extend ${\DRconj}$ to
 $\F^{-(k+1)}(Y\cap X)$ in such a
 way that the extension is continuous 
  along rays. Indeed, for every
 $\s$, we can choose a branch $L$ 
  of $\Ek^{-1}$ on the ray such that
 $L(\DRconj(\F(\s,t)))=\DRconj(\s,t)$ whenever 
   $(\s,t)\in \F^{-k}(Y\cap X)$.
 This extension is also continuous in 
  both variables because the branch $L$
 varies continuously.
\end{proof}

As a direct consequence of Corollary \ref{cor:globalcorrespondence}, 
 we obtain the classification
 theorem from \cite{expescaping}.
\begin{cor}[Classification of Escaping Points 
            \protect{\cite[Corollary 6.9]{expescaping}}] \label{cor:classification}
 Let $\Ek$ be an exponential map. For every escaping point $z\in
 I(\Ek)$, exactly one of the following holds:
 \begin{itemize}
  \item There exists a unique $x\in X$ such that
         that $z=\DRconj(\s,t)$, or
  \item the singular value $\kappa$ escapes; there exist $\s$ and
         $t_0 > \ts$ with $\kappa=\DRconj(\s,t_0)$, and there is $n\geq 1$
         such that $\Ek^n(z)=\DRconj(\s,t)$ with $\ts \leq t \leq t_0$. \qedd
 \end{itemize}
\end{cor}

\begin{proof}[Proof of
  Theorem \ref{thm:boettcheranalog}]
  Set $K:= 
   \max(|\kappa_1|,|\kappa_2|,
		      2\pi+6)+1$,
   and
   define $R:= Q(K)+1$. 
   For every $\kappa\in\C$
    with $|\kappa|<R$,
    let $\DRconj^{\kappa}:Y_{Q(K)}\to\C$ denote
    the map from Theorem
    \ref{thm:conjugacyconvergence}. 

  If $A$ is as in the statement of
   the theorem, then by
   Theorem \ref{thm:surjectivity},
    $A\subset \DRconj^{\kappa_1}(Y_{Q(K)})$.
    Thus the map
   \[ \Phi^{\kappa}:A\to\C, z\mapsto
           \DRconj^{\kappa}\bigl(
            (\DRconj^{\kappa_1})^{-1}(z)\bigr) \]
    conjugates $E_{\kappa_1}$ on
    $A$ to $E_{\kappa}$ on
    $\Phi^{\kappa}(A)$. Since $\Phi$ is
    a holomorphic motion of
    $A$ (compare \cite{mss}),
    $\Phi^{\kappa}$ extends to
    a quasiconformal homeomorphism
    $\C\to\C$.

   If $\kappa_1$ and $\kappa_2$
    are nonescaping parameters, then
    $\Phi^{\kappa}|_{A\cap I(f)}$
    extends to the bijection
    $\DRconj^{\kappa_2}\circ \DRconj^{\kappa_1}$
    from Corollary 
    \ref{cor:globalcorrespondence}
    (which is continuous along rays,
     but not continuous in general).
  \end{proof}

When considering individual rays, it is often
 cumbersome to take into account the starting potential
 $t_{\s}$. For convenience, we make the
 following definition. 
\begin{defn}[Dynamic Rays]
 Let $\Ek$ be an exponential map and let $\s\in\Sequ_0$.
 We define a curve $g_{\s}$ --- the
 \emph{dynamic ray at address $\s$} --- by
  \[ g_{\s}(t) = {\DRconj}(\s,t+t_{\s}). \]

 If $\gs$ is not defined for all $t>0$ (i.e., if there exists
  $t_0>\ts$ such that ${\DRconj}(\F^n(\s,t_0))=\kappa$), then we
  call $\gs$ a ``broken ray''.
  We say that an unbroken ray $g_{\s}$ \emph{lands} 
  at a point $z_0$ if $\lim_{t\to 0}\gs(t)=z_0$. Similarly, we say that
  $\gs(t)$ \emph{has an escaping endpoint} if $\gs(0)$ is defined and
  escaping; i.e. if $(\s,\ts)\in X$.
\end{defn}

\begin{figure}
 \subfigure[A beginning piece of the curve $\widetilde{g_{\per{10}}}$
            obtained when approximating $g_{\per{10}}$
            from above. $\widetilde{g_{\per{10}}}$ first traverses the ray
            $g_{\per{10}}$ (the bottom curve in the picture), from right
            to left, then the upper curve (a preimage of the ray piece
            connecting $\Ek(\kappa)$ to $-\infty$) from left to right,
            followed by further preimages of this piece.]{%
            \includegraphics[width=.89\textwidth]{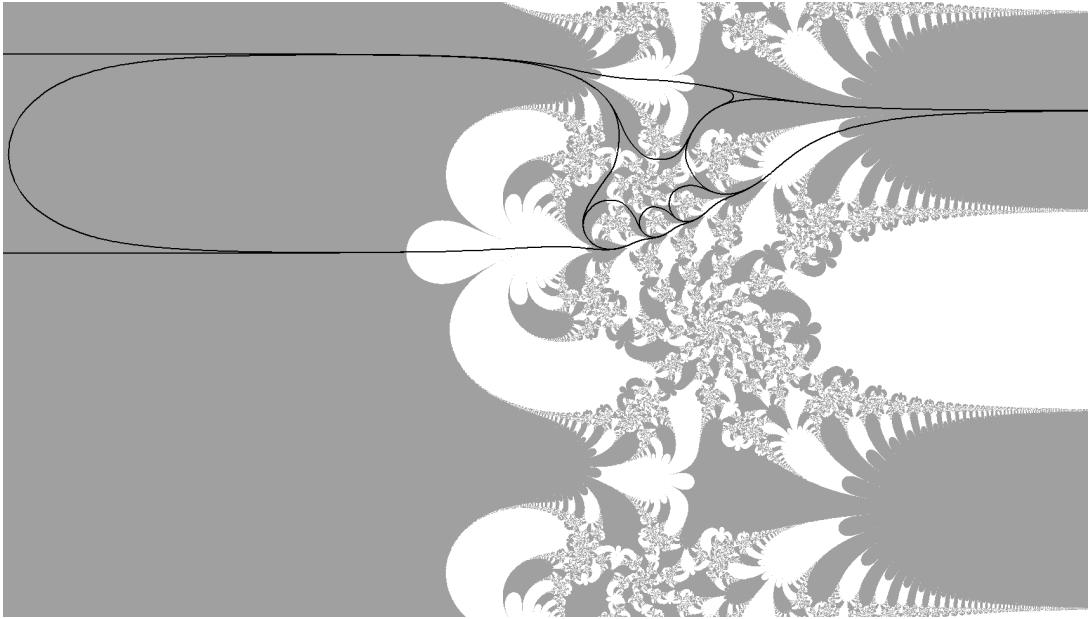}}

 \vfill

 \noindent
 \subfigure[The analogous situation when approximating $g_{\per{10}}$ from below.]{%
            \includegraphics[width=.89\textwidth]{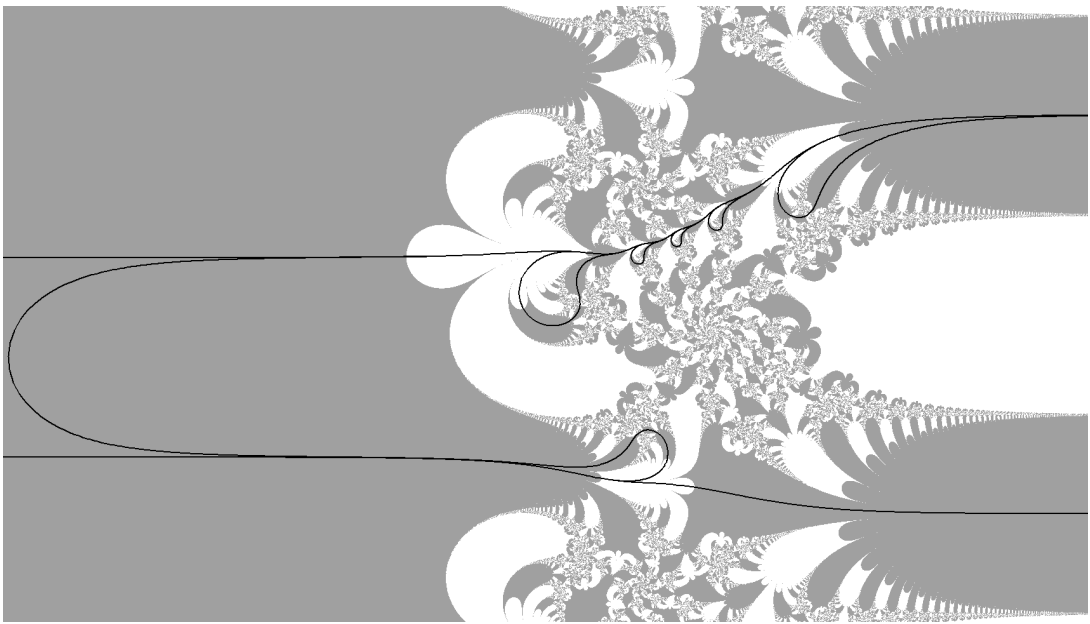}}
 \caption{An illustration of the remark about broken rays after Lemma
            \ref{lem:raycontinuity}. Here $\kappa\in \gs$, where $\s$ is the
            periodic address $\s=\per{01}$. 
          \label{fig:brokenrays}}
\end{figure}

\begin{lem}[Convergence of Rays] \label{lem:raycontinuity}
 Let $\Ek$ be an exponential map. 
 Suppose that $\s^n$ 
 is a sequence of external addresses converging to
 a sequence $\s^0\in\Sequ_0$ 
 such that also $t_{\s^n}\to t_{\s^0}$, and let $t_0>0$
 such that $g_{\s^0}(t)$ is defined for all $t\geq t_0$. Then
 \[ g_{\s^n}|_{[t_0,\infty)} \to g_{\s^0}|_{[t_0,\infty)} \]
 uniformly.
\end{lem}
\begin{proof} 
 Let $Q=Q(|\kappa|)$ as before. 
 There exists $k$ such that $\{(\s,t):t - \ts \geq t_0\}\subset
 \F^{-k}(Y_Q)$. The claim then follows from Corollary
 \ref{cor:globalcorrespondence}.
\end{proof}
\begin{remark}
  In the case where $g_{\s^0}$ is broken, we can say the
 following (with the
 same proof). Suppose that $\s^n> \s^0$ (or $< \s^0$) for all
 $n$. Then, under the assumptions of Lemma \ref{lem:raycontinuity}
 the rays $g_{\s^n}$ converge locally uniformly (in the spherical
 metric on the Riemann Sphere $\Ch$) to a curve
 $\widetilde{g_{\s^0}}:(0,\infty)\to \Ch$. This curve has
 $\widetilde{g_{\s^0}}(t)=\infty$ if and only if
 ${\DRconj}(\F^n(\s^0,t+t_{\s^0}))=\kappa$ for some $n\geq 1$ and 
 coincides with
 $g_{\s^0}$ where the latter is defined. If the ray which
 contains $\kappa$ is periodic, then $\infty$ is assumed infinitely many
 times on this curve. In this case, the curve
 accumulates everywhere on itself. (See Figure \ref{fig:brokenrays}.)
 A Theorem of Curry \cite{curry} can then be used to show that
 the accumulation set of $\widetilde{g_{\s^0}}$ in
 $\C$ can be compactified to an
 indecomposable continuum. This was previously done for $\kappa\in (-1,\infty)$ in
 \cite{devaneyknaster} and for certain other parameters in \cite{moreno}.
\end{remark}

In the remainder of the article, we sometimes write ${\DRconj}^{\kappa}$ or
 $g^{\kappa}_{\s}$ for the objects constructed previously when the
 parameter is not fixed in the context.

There is an interesting
corollary of Theorem \ref{thm:conjugacyconvergence}.
Note that $Y_Q$ contains many points of
$\Xb\setminus X$; in particular endpoints of periodic addresses.
Which of these addresses lie in $Y_Q$ depends on $Q$
(and thus on the parameter); 
however, we can use this fact to give an elementary bound on
those parameters for which we know that these rays cross sector
boundaries or are not defined. This result is used in
\cite{landing2new} and \cite{boundary} to bound \emph{parameter
rays} (see Section \ref{sec:parameterspace}) and \emph{wakes of
hyperbolic components}.

\begin{cor}[Bound on Parameter Rays] \label{cor:parameterraybound} 
 Let $\kappa\in\C$, $Q := Q(|\kappa|)$ as in
  (\ref{eqn:Q}) and suppose that $(\s,t)\in\Xb\setminus Y_Q$. If
  the number $t_0 := \inf_{j\geq 0} T(\F^j(\s,t))$ satisfies
  $t_0\geq \pi+2$, then
  $|\kappa| > \frac{1}{5}\exp(t_0)$. 

 In particular, suppose that $\s$ is periodic
  of period $n$ and $M := \max |s_k|\geq F^{n}(\pi+2)/2\pi$. If
  $\kappa\in\C$ with $\kappa\in g_{\s}^{\kappa}$, then
  $|\kappa| \geq \frac{1}{5}F^{-(n-1)}(2\pi M)$. 

 Similarly, suppose that $\s^1,\s^2\in\Sequ_0$ 
  are external addresses for
  which there
  is $n\in\N$ and $M\geq F^{n}(\pi+2)/2\pi$ such that 
 \[ \max_{k+1\leq m \leq k+n} |\s^j_m| \geq M \]
 for all $k\geq 0 $ and $j\in\{1,2\}$. If $\kappa$ is a parameter such
 that $g^{\kappa}_{\s^1}$ and $g^{\kappa}_{\s^2}$ land together, then
 $|\kappa| \geq \frac{1}{5}F^{-(n-1)}(2\pi M)$.
\end{cor}
\begin{proof}
 The first claim is an immediate consequence of 
  the definition of $Q$ and $Y_Q$. Note that,
  if $\s$ is an address such that
  $\kappa\in g_{\sigma(\s)}$ or such that
  $g_{\s}$ lands at a point which does not have 
  external address $\s$, then $(\s,\ts)\notin Y_Q$. 

 Furthermore, among all addresses 
 $\s$ one of whose entries $s_2,\dots,s_{n+1}$ 
  is of size at
  least $M$, the value of $\ts$ is minimized by the address
   $\s = 0 0 \dots 0 M \per{0}$
  (where the first block of $0$s consists of $n$ entries).
  For this $\s$, we have
  $\ts = F^{-n}(2\pi M)$. Thus the second and third statements follow from the
  first.
\end{proof}

\section{Limiting Behavior of Dynamic Rays}
 \label{sec:limitsofrays}

 For completeness, this section collects some results
  from \cite{expescaping} and \cite{markuslassedierk} on the
  limiting behavior of dynamic rays. First we state and prove
  two lemmas which were implicitly contained in the proof of
  \cite[Corollary 6.9]{expescaping} and imply that
  a dynamic ray at a slow external address cannot land at
  an escaping point. From this, we deduce, as
  first outlined in \cite{markuslassedierk}, that every
  dynamic ray is a path-connected
  component of $I(\Ek)$. 

\begin{lem}[Limit Set of Ray] \label{lem:limitset}
 Let $g:(0,\infty)\to I(\Ek)$ be an unbroken dynamic ray, and let
   \[ L := \bigcap_{t>0} \cl{g\bigl((0,t)\bigr)} \]
  denote the limit set of $g$. If there is $(\s,t_0)$ such that
  $g_{\s}(t_0)\in L$, then $\gs(t)\in L$ whenever $t\leq t_0$ is
  such that $\gs(t)$ is defined.
\end{lem}
\begin{remark}
 If the ray $\gs$ is broken, we could replace $\gs$ in the last
  statement by the curve $\wt{\gs}$ from the remark after Lemma
  \ref{lem:raycontinuity}.
\end{remark}
\begin{proof}
  Let us define addresses $\s^{n\pm}$ by
   \[ s^{n\pm}_{k} = \begin{cases}
                        s_k \pm 1 & k = n \\
                        s_k & \text{otherwise}.
                      \end{cases}
   \]
       One sees easily that
        $t_{\s^{n\pm}}\to t_{\s}$ for $n\to\infty$
    (compare e.g.\ Lemma \ref{lem:tsstar}).
   Now pick some $t$, $0<t<t^0$. Then, by Lemma \ref{lem:raycontinuity},
    \[g_{\s^{n\pm}}\bigl([t,\infty) \bigr) \to
      g_{\s}\bigl([t,\infty) \bigr) \]
   uniformly. Therefore any curve which does not intersect the
    $g_{\s^{n\pm}}$ and accumulates at $g_{\s}(t_0)$ must also
    accumulate at
    $g_{s}(t)$.
\end{proof}

\begin{lem}[Addresses of rays landing together]
   \label{lem:rayslandingtogether}
 Suppose that for some $\s,\s'\in\Sequ_0$ the dynamic rays
 $g_{\s}$ and $g_{\s'}$ land at the same point. Then $|s_k - s_k'|\leq
 1$ for all $k$.
\end{lem}
\begin{proof}
 Assume, by contradiction, that $|s_k- s_k'|>1$ for some $k$; by
 passing to a forward iterate
 if necessary we can assume that $k=1$. 
 Let $S$ denote the union of $g_{\s}$, $g_{\s'}$ and their
 common landing point $z_0$, which is a Jordan arc tending to
 $\infty$ in both directions. Note that $\Ek$ is injective on
 $S$. Indeed, $\Ek$ is injective on every ray, and it is injective on
 $\gs\cup g_{\s'}$ unless $\s$ and $\s'$ differ only in their first
 entries. However, in that case $g_{\s'}$ would be a translate of
 $\gs$ by a multiple of $2\pi i$, which means that $\gs$ and $g_{\s'}$
 cannot land together. Finally, $\Ek$ is injective on $S$ as otherwise
 $z_0$ would lie on a dynamic ray, which contradicts Lemma
 \ref{lem:limitset}. 

 On the other hand, if $|s_1 - s_1'|>1$, then the two ends of $S$ tend
 to $\infty$ with a difference of more than $2\pi$ in their imaginary
 parts, which implies that $S\cap (S+2\pi i ) \neq \emptyset$. Thus
 $\Ek$ is not injective on $S$, a contradiction.
\end{proof}

\begin{cor}[No Landing at Escaping Points] \label{cor:nolandingatescapingpoints}
 Suppose that $\s\in\Sequ_0$ is a slow external address (i.e.\
  $(\s,\ts)\notin X$). Then $g_{\s}$ does not land at an
  escaping point. 
\end{cor}
\begin{proof} By Lemma \ref{lem:limitset}, 
  $g_{\s}$ could only land at the escaping endpoint
 $g_{\s'}(t_{\s'})$ for some fast address $\s'$.
 By Lemma \ref{lem:rayslandingtogether}, $|s_k - s_k'|\leq 1$ for all
 $k$. It easily follows that $\s$ would also have to be a fast
 address, which contradicts our assumption. (Compare 
 Corollary \ref{cor:properties}.)
\end{proof}

To infer that the path-connected components of $I(\Ek)$ are given
 by dynamic rays, we require the following topological fact.

 \begin{prop}[Path Components {\cite[Proposition 4.2]{markuslassedierk}}]
     \label{prop:pathcomponentsarecurves}
  Let $I$ be a Hausdorff topological space. Let $\Gamma$ be
   a partition of $I$ into path-connected subsets such that
   no union of two different elements of $I$ is path-connected.

  Suppose that $I$ can be written as a countable union of closed subsets
   $I_k$ such that
  \begin{enumerate}
   \item $I_{k}\subset I_{k+1}$
   \item every path-connected component of $I_k$ is
     contained in some element of $\Gamma$, 
      \label{item:pathcomponentsofIk}
   \item every element of $\Gamma$ contains at most 
     one path-connected component of $I_k$, and \label{item:onlyonecomponent}
   \item \label{item:complicated}
     if $c\subset I$ is a simple closed
     curve which is completely contained in some element
     of $\Gamma$, then either $c\subset I_k$ or $c\cap I_k = \emptyset$.
  \end{enumerate}

  Then $\Gamma$ is the set of path-connected components of $I$. \qedd
 \end{prop}

 \begin{cor}[Path Components of $I(\Ek)$ 
               {\cite[Corollary 4.3]{markuslassedierk}}]
  Let $\kappa\in\C$. Then every path connected component of
  $I(\Ek)$ is
  \begin{enumerate}
   \item a dynamic ray, \label{item:slowray}
   \item a dynamic ray together with its escaping endpoint, or 
     \label{item:fastray}
   \item (if $\kappa\in I(\Ek)$) an iterated preimage component of
       the dynamic ray containing the singular value. \label{item:brokenray}
  \end{enumerate}
 \end{cor}
 \begin{proof}
   Set $I:=I(\Ek)$, and let $\Gamma$ be the set of curves of types
   (\ref{item:slowray}) to (\ref{item:brokenray}). By Corollary
   \ref{cor:nolandingatescapingpoints}, no union of two different
   elements of $\Gamma$ is path-connected, and no element of $\Gamma$
   contains a simple closed curve. 

   We now set $I_0 := {\DRconj}(Y_Q\cap X)$ and
   $I_{j} := \Ek^{-j}(I_0)$. It is easy to see that these sets
   satisfy the hypotheses of Proposition \ref{prop:pathcomponentsarecurves},
   and the claim follows. 
 \end{proof}

Little else is known about the possible limiting behavior of dynamic
 rays. For polynomials, dynamic rays cannot accumulate at escaping
 points. This is not the case for exponential maps: in
 \cite{indecomposable} it was shown that, for $\kappa\in (-1,\infty)$,
 there exists a ray which accumulates on itself. This was also shown for
 a larger class of exponential maps in \cite{nonlanding};
 see also \cite[Section 3.8]{thesis}. 
 On the other hand, it is now known that every (pre-)periodic dynamic
 ray of an exponential map lands \cite{landing2new}.

\begin{figure}
 \fbox{%
 \includegraphics[width=.3\textwidth]{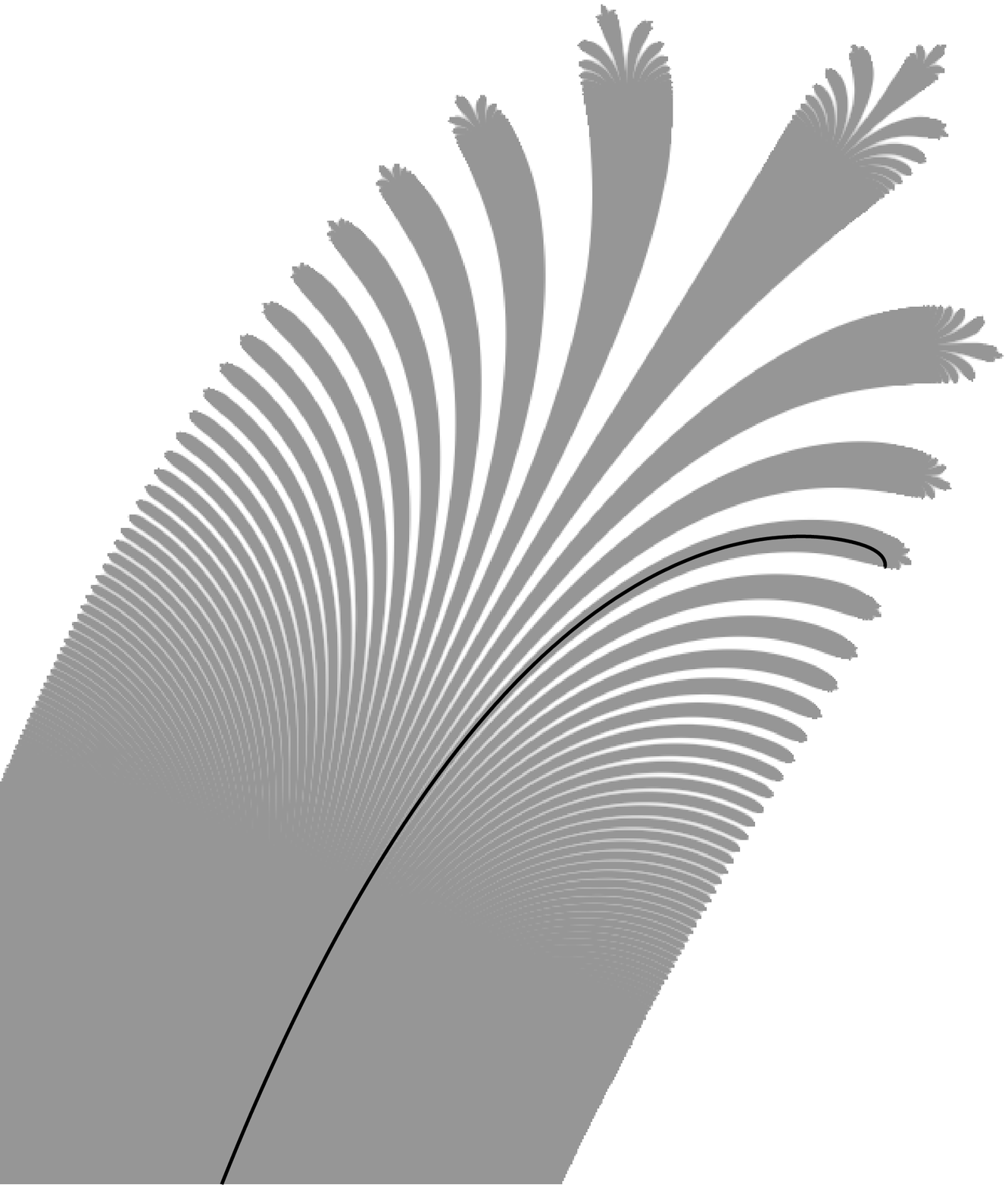}}\hfill
 \fbox{%
 \includegraphics[width=.3\textwidth]{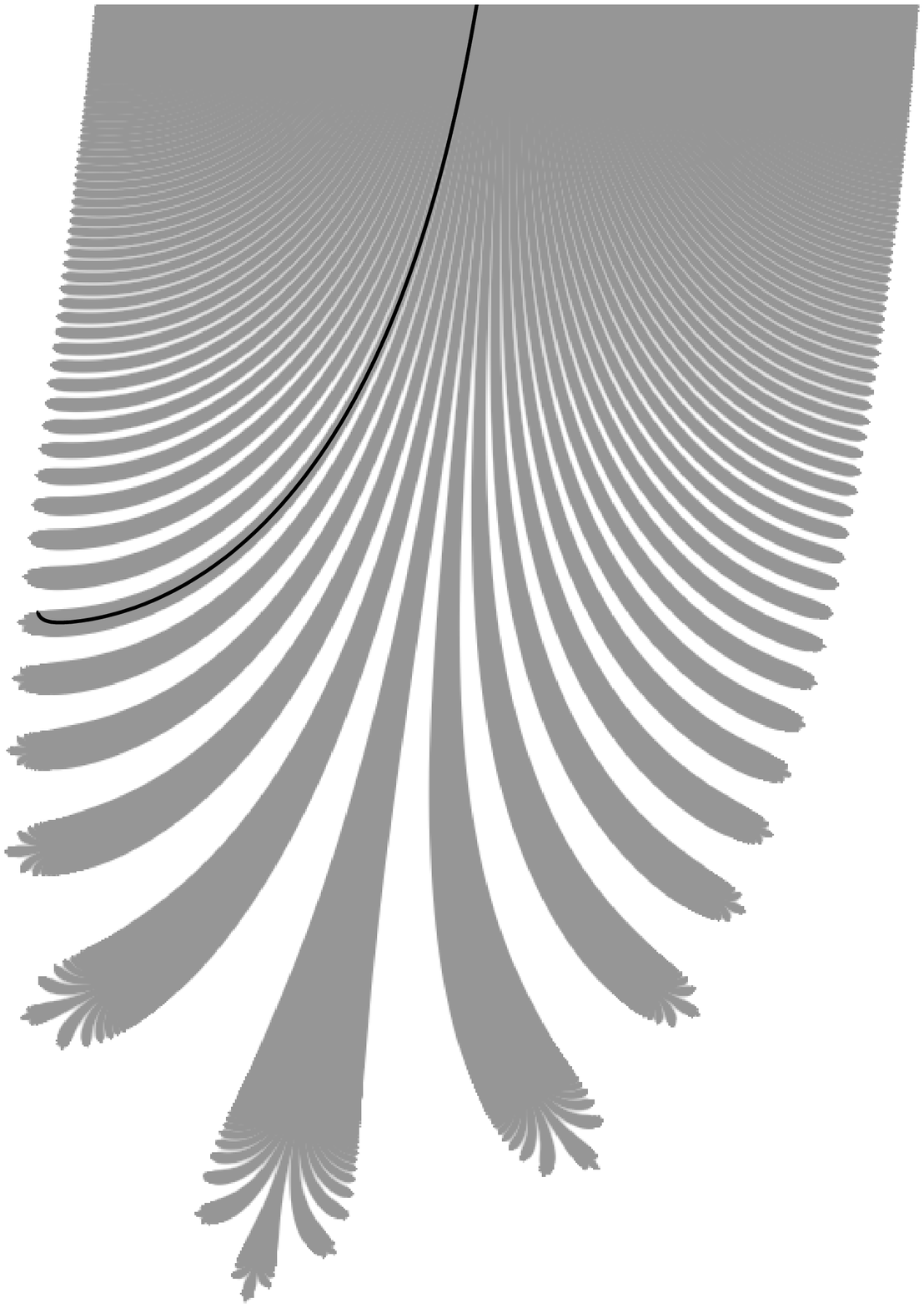}}\hfill
 \fbox{%
 \includegraphics[width=.3\textwidth]{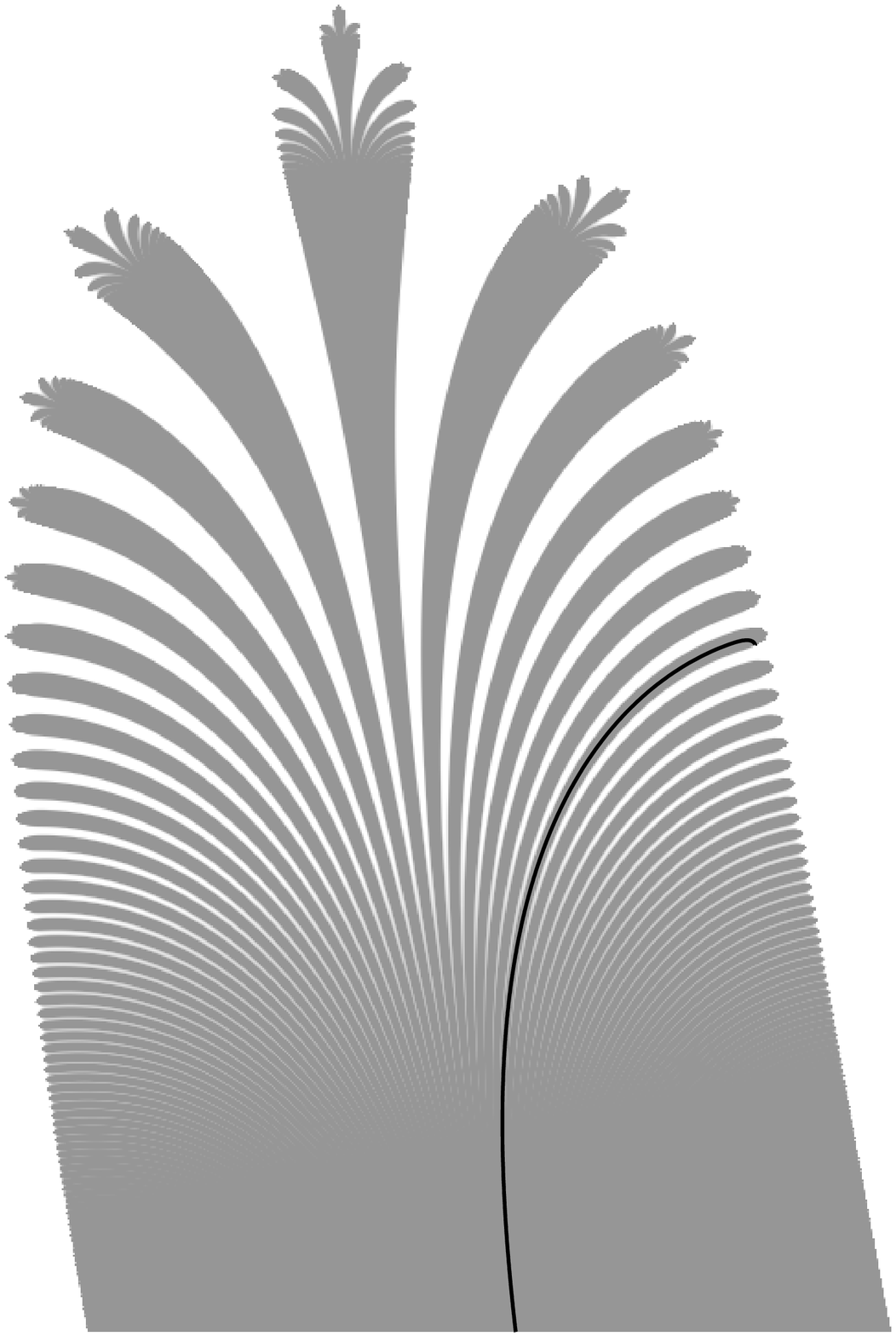}}
 \caption[A Ray Not Differentiable in its Endpoint]{%
  These pictures are intended to illustrate Theorem
 \ref{thm:endpointdifferentiability}. Shown is the ray at address
 $01234\dots$ for the parameter $\kappa=-2$ in three successive
 magnifications. The images demonstrate that the ray indeed spirals around
 its escaping endpoint.}
\end{figure}

\section{Differentiability of Rays} \label{sec:differentiability}

Viana \cite{viana} proved (using a different parametrization) 
that the rays $\gs$ are
$C^{\infty}$. His arguments also apply to the parametrization of the
curves given by our construction. (Compare also the proof of Theorem
\ref{thm:endpointdifferentiability} below.)

\begin{thm}[Rays are Differentiable \cite{viana}]
 Let $\s\in\Sequ_0$. Then
 $\gs:(0,\infty)\to\C$ is a $C^{\infty}$ function. \qedd
\end{thm}

A proof of the differentiability of rays
can also be found in \cite{markusdierk},
where this was carried out to obtain specific estimates on the first
and second derivatives. However, previously there was no information
about which rays with escaping endpoints are also
differentiable in these endpoints, and whether this may depend on
the parameter. Using the results of
Section \ref{sec:classification}, we can answer this question.

\begin{thm}[Differentiability of Rays in Endpoints]
\label{thm:endpointdifferentiability} 
 Let $\s\in\Sequ_0$ be a fast external address,
  and let $\kappa$ be a parameter for which the ray
  $\gs$ is unbroken. Then the curve $\gs\bigl([0,\infty)\bigr)$ 
  is continuously differentiable in $\gs(0)$
  if and only if the series
  \begin{equation}
    \sum_{j=0}^{\infty} \frac{2\pi s_{j+1}}{T(\F^j(\s,\ts))}
    \label{eqn:argumentsum}
  \end{equation}
  converges.
\end{thm}
\begin{remark}
  By the formulation ``the curve is continuously differentiable in $\gs(0)$''
  we mean that it is continuously differentiable under a suitable
  parametrization (e.g., by arclength), \emph{not} that the function
  $\gs$ itself is necessarily differentiable in $0$. 
  If the convergence of
  the sum is absolute, then one can show that the function $\gs$ itself is
  differentiable in $0$.
\end{remark}
\begin{proof}
  Let $Q$ be as in 
  (\ref{eqn:Q}). Then 
  all
  ${\DRconj}_k$ are defined on the set $Y_Q$ and converge uniformly to the
  function ${\DRconj}$ there; it is clearly sufficient to prove the theorem
  for addresses for which $(\s,\ts)\in Y_Q\cap X$.

  By the definition of the functions ${\DRconj}_k$, their $t$-derivatives
  in any point $(\s,t)\in Y_Q$ are given by
  \begin{align*}
    \frac{\partial {\DRconj}_k}{\partial t} (\s,t) &=
       \frac{1}{{\DRconj}_{k-1}(\F(\s,t)) - \kappa} \cdot
        \frac{\partial {\DRconj}_{k-1}}{\partial t}(\F(\s,t)) \cdot \exp(t) =
         \dots \\
      &= \prod_{j=1}^k
   \frac{\exp(T(\F^{j-1}(\s,t)))}{{\DRconj}_{k-j}(\F^j(\s,t)) - \kappa} \\
      &= 
   \left( \prod_{j=1}^k\frac{\exp(T(\F^{j-1}(\s,t)))}{{\DRconj}(\F^j(\s,t)) - \kappa}
     \right) \cdot
       \prod_{j=1}^k\left( 1 + 
                     \frac{{\DRconj}(\F^j(\s,t))- {\DRconj}_{k-j}(\F^j(\s,t))}{%
                           {\DRconj}_{k-j}(\F^j(\s,t)) - \kappa} \right).
  \end{align*}

Recall from the proof of Theorem \ref{thm:conjugacyconvergence} that
  \[ \bigl|{\DRconj}(\F^j(\s,t))- {\DRconj}_{k-j}(\F^j(\s,t))\bigr| \leq
     2^{-(k-j)} \cdot (2\pi+4), \]
 so the second product converges uniformly for $t\geq \ts$. It is not
  difficult to see that the first product converges locally uniformly
  (and is nonzero) 
  for $t > \ts$ (see e.g. \cite{markusdierk}). Note that this proves that
  the ray
  without the endpoint is $C^1$.

 The ray is continuously differentiable in its endpoint if and
 only if $\arg \left( \frac{\partial G}{\partial t} (\s,t)\right)$ 
 has a limit as $t\to
 \ts$. 
 The above argument shows that this is equivalent to the
 question whether
 the function 
  \[ \Theta(t) := \sum_{j=1}^{\infty} \arg\bigl( {\DRconj}(\F^j(\s,t))-\kappa \bigr) \]
 has a limit for $t\to \ts$. Let us set
 Let us set $\arg(k,t):=\arg(\DRconj(\F^k(\s,t))-\kappa)\in
                                (-\pi,\pi)$ for all $k\geq 1$ and
    $t\geq \ts$. 
 \begin{claim} 
  The limit $\lim_{t\to\ts} \Theta(t)$ exists if and only if
   the series $\Theta(\ts)$ is convergent.
 \end{claim}
 To prove this, choose some number $m\geq 3$ (to be fixed below) and
  define, for $n$ large enough, $t_n>\ts$ to be
  the unique number for which
  \[ T(\F^n(\s,t_n)) = T(\F^n(\s,\ts)) + \log m . \]
                                                                     
 We will prove the claim by comparing the summands of $\Theta(t_n)$
    with those of $\Theta(t_{\s})$.
  Note that $|{\DRconj}(\F^n(\s,t_n)) - {\DRconj}(\F^n(\s,\ts))| \leq
     K := 2\pi + 4 +
    \log m$. It follows again by contraction that, for $k\leq n$,
  \[ |{\DRconj}(\F^k(\s,t_n)) - {\DRconj}(\F^k(\s,\ts))| \leq 2^{-(n-k)}\cdot K. \]
  Thus
  \begin{align*}
    \biggl|\sum_{k-1}^n \arg(k,\ts) - \sum_{k=1}^{n}\arg(k,t_n)\biggr| 
     \leq 
    \pi K &\sum_{k=1}^n \frac{2^{-(n-k)}}{|{\DRconj}(\F^k(\s,\ts)) -\kappa|},
  \end{align*}
  which is easily seen to converge to $0$ as $n\to\infty$. 

 Also observe that, for $k\geq n+1$, 
  \begin{align*}
    T(\F^k(\s,t_n)) &\geq
       F^{k-n-1}\bigl(T(\F^{n+1}(\s,t_n)) -
                      T(\F^{n+1}(\s,\ts))\bigr) \\
     =& F^{k-n-1}\bigl((m-1)\cdot \exp(T(\F^n(\s,\ts)))\bigr)
     \geq F^{k-n-1}\bigl((m-1)\cdot F(T(\F^{n}(\s,\ts)))\bigr)
  \end{align*}
 and
  \[ 2\pi s_{k+1} \leq F^{k-n-1}\bigl(F(T(\F^{n}(\s,\ts)))\bigr). \]
 It easily follows that
  \[ \sum_{k=n+2}^{\infty}\bigl|
     \arg(k,t_n)\bigr|  \to 0 \]
  as $n\to\infty$. Similarly, for large enough $n$,
  the value
  $\bigl|\arg(n+1,t_n)\bigr|$
  is no larger than
  $\frac{2}{m-1}+\eps$. If $\arg(n+1,\ts)$
  tends to $0$ as $n\to\infty$, then 
  $\arg(n+1,t_n)$ also does.

 Now let us first consider the case that $\arg(n,\ts)$
  does
  not converge to $0$ (and thus the sum $\Theta(\ts)$ is
  divergent). So let $\delta>0$ and let $n_k$ be a subsequence for which
  $\arg(n_k,\ts)\geq \delta$.
  If $m$ was chosen to be $1+\frac{5}{\delta}$, then it follows
    from the above considerations that
   \[ |\Theta(t_{n_k-1}) - \Theta(t_{n_k})| \geq
       \delta - \frac{4}{m-1} + o(1) > \frac{\delta}{5} + o(1) \]
   (as $k\to\infty$). 
  In particular, the sequence $\Theta(t_n)$ 
  does not have a limit for $n\to\infty$. This proves the claim in this
  case.

 So we can now suppose that $\arg {\DRconj}(\F^j(\s,\ts))\to 0$. Then, by our
  observations, 
  \[ \bigl| \Theta(t_n) - 
     \sum_{k=1}^n \arg(k,\ts) \bigr| \to
     0. \]
 Thus in particular the sequence $\Theta(t_n)$ has a limit if and only 
  if the sum $\Theta(\ts)$ is convergent. It remains to show that this 
  implies that $\Theta$ has a limit as $t\to \ts$. However, it is easy
  to show that
    \[ \sup_{t\in[t_n,t_{n+1}]} |\Theta(t_n) - \Theta(t)| \to 0 \]
  as $n\to \infty$. Indeed, by the above observations,
    \[ \sum_{k=1}^{n-1} \left| \arg(k,t) -
                               \arg(k,t_n) \right| \]
  is small, as is
    \[ \sum_{k=n+2}^{\infty} |\arg( {\DRconj}(\F^{k}(\s,t)) - \kappa )|. \]
  The two entries that remain to be dealt with tend to $0$ because
  $\arg {\DRconj}(\F^j(\s,\ts))$ does. 
  This proves the claim in the second case. 

 To conclude the proof, we need to show that the convergence 
 of the sum $\Theta(\ts)$ is equivalent to the convergence of the sum
 (\ref{eqn:argumentsum}) in the statement of the theorem. 
 It is clear that the terms of
 (\ref{eqn:argumentsum}) converge to $0$ if and only if those of 
 $\Theta(\ts)$ do. So we
 can suppose that $\arg {\DRconj}(\F^k(\s,\ts))\to 0$. It is easy to see that
 then there exists $x>0$ such that $|{\DRconj}(\F^k(\s,\ts))|\geq F^k(x)$ for
 all large enough $k$. (Compare Corollary \ref{cor:properties}.)
 Because of this and since 
 $|{\DRconj}(\F^k(\s,\ts)) - Z(\F^k(\s,\ts))|$ is bounded by $2 + \pi$, we have
 \[\left| \sum_k^{k_0}  \arg ({\DRconj}(\F^k(\s,\ts))-\kappa) - 
           \sum_k^{k_0} \arg Z(\F^k(\s,\ts)) \right| \leq
        \pi\cdot(2+\pi+|\kappa|)\cdot \sum_k^{k_0}
 \frac{1}{F^k(x)}. \]
 The last sum is clearly absolutely convergent, and so the convergence of the
 sum $\Theta(\ts)$ and that of
  \[ \sum_k^{\infty} \arg Z(\F^k(\s,\ts)) \]
 are equivalent. Similarly, one sees that the convergence of this last
 sum and the sum (\ref{eqn:argumentsum}) are equivalent. 
\end{proof}


\section{Speed of Escape} \label{sec:speedofescape}

 Our construction
  provides information on the speed of
  escaping endpoints of dynamic rays, about which previously
  little was known. To deal with such questions, it is often useful 
  to introduce a quantity $\ts^*$ which is closely related
  to $\ts$, but can be computed more
  easily.

\begin{deflem}[Growth of Escaping Endpoints] \label{lem:tsstar}
 Let $\s$ be an external address and  define
    \[ \ts^* := \sup_{k\geq 1} F^{-k}(2\pi|s_{k+1}|). \]
  Then $\s$ is exponentially bounded if and only if
    $\ts^*<\infty$, in which case $\ts^* \leq \ts \leq \ts^*+1$.
\end{deflem}
\begin{proof}
 Suppose first that $\s$ is exponentially bounded; i.e.\ $\ts<\infty$. 
  By definition of $\ts$, we have
  $F(\ts)\geq 2\pi|s_2|$ for all external addresses $\s$. 
  Since $t_{\sigma(\s)}\leq F(t_{\s})$, it follows inductively
  that $F^k(\ts)\geq 2\pi|s_{k+1}|$. This proves that
  $\ts^*\leq \ts < \infty$.

  Now suppose that $\ts^*<\infty$. Then
   \[ T\bigl(\F(\s,\ts^*+1)\bigr) =
       F(\ts^*+1)-2\pi|s_2| \geq
       2F(\ts^*)+1-2\pi|s_2| \geq 
       F(\ts^*)+1 \geq t_{\sigma(\s)}^*+1 \]
   (where we used the facts that $F(t+1)\geq 2F(t)+1$ and
    $2\pi|s_2|\leq F(t_{\s}^*)$). 
   It follows by induction that $(\s,\ts^*+1)\in\Xb$; i.e.\
   $\ts \leq \ts^*<\infty$. 
\end{proof}

 As a first application, we recover the characterization of
  exponentially bounded, slow and fast addresses given
  in \cite{expescaping}. Also, we obtain a new
  description of addresses with
  \emph{positive minimal potential} in the sense of
  \cite{expescaping}. These are addresses $\s\in\Sequ_0$
  for which there exists $x>0$ with $2\pi|s_n|\geq F^{n-1}(x)$ 
  for infinitely many $n$.

\begin{cor}[Properties of External Addresses] \label{cor:properties}
 Let $\s$ be an external address. Then
 \begin{itemize}
  \item $\s$ is exponentially bounded if and only if
         there exists $x\geq 0$ with $2\pi|s_k|\leq F^{k-1}(x)$ for
         all $k\geq 1$,
  \item $\s$ is slow if and only if there exist $x\geq 0$ and infinitely
         many $n\geq 0$ which satisfy $2\pi|s_{n+k}|\leq F^{k-1}(x)$ for all 
         $k\geq 1$, and
  \item $\s$ has positive minimal potential
          if and only if
          there exists $x>0$ with the property that
          $T(\F^{j}(\s,t_{\s}))=t_{\sigma^j(\s)}>F^j(x)$ for all $j$.
 \end{itemize}
\end{cor}
\begin{proof} By the
  previous lemma, $\s$ is exponentially bounded if and only if 
  $\ts^* <\infty$,
  which is clearly equivalent to the stated condition. 
  The other claims follow in a similar fashion. 
\end{proof}

 We now prove Theorem
  \ref{thm:arbitrarilyslowescape}, restated here for convenience.

 \begin{thm}[Arbitrary Escape Speed]
  Let $\kappa\in\C$. 
  Suppose that $r_n$ is a sequence of positive real numbers such that
  $r_n\to \infty$ and $r_{n+1}\leq \exp(r_n) + c$ for some $c>0$. Then
  there is an escaping endpoint
  $z\in I(\Ek)$ and some $n_0\in\N$ such that, for
  $n\geq n_0$,
   $|\re(\Ek^{n-1} (z)) - r_n | \leq 2+2\pi$.
 \end{thm}
 \begin{proof}
  By changing the first few entries of $(r_n)$, if necessary,
   we may assume without loss of generality that 
   $F(r_n+1) > r_{n+1}+1$. We define an external address $\s$ by
   $s_{n+1} := \lfloor F(r_{n})/2\pi\rfloor$. Then
   \[ t_{\sigma^{n-1}(\s)}\geq F^{-1}(2\pi|s_{n+1}|) \geq 
       r_n - 2\pi \]
   for all $n$. Furthermore
   \[ F^{-k}(2\pi|s_{n+k}|) \leq
      F^{-(k-1)}(r_{n+k-1}) \leq
      r_n+1 \]
   for $k\geq 1$. 
   By Lemma \ref{lem:tsstar}, we thus have
     \[ r_n - 2\pi \leq 
         t_{\sigma^{n-1}(\s)} \leq
         t_{\sigma^{n-1}(\s)}^* + 1 \leq r_n + 2 \]
   for all $n$. The claim now follows by Lemma
    \ref{lem:boundeddistance} and Theorem
   \ref{thm:conjugacyconvergence}: we pick $n_0$
   sufficiently large so that $\F^{n_0}(\s,\ts)\in Y_Q$ and
   $z_{n_0} := \DRconj(\F^{n_0}(\s,\ts))\neq\kappa$ 
   and choose $z$ to be any element of
   $\Ek^{-n_0}(z_{n_0})$. 
 \end{proof}

 Frequently, one considers escaping points which eventually
  escape in a sector (i.e., $\im \Ek^n(z)\leq C\re\Ek^n(z)$)
  or, more generally, a parabola (i.e., $\im \Ek^n(z)\leq C(\re\Ek^n(z))^K$).
  For example, McMullen \cite{hausdorffmcmullen}
  showed that the set of escaping points satisfying a sector condition
  has Hausdorff dimension two. On the other hand,
  Karpinska \cite{karpinska}
  proved that the Hausdorff dimension of the set of points
  satisfying a parabola condition of exponent $K<1$ is at most
  $1+K$. Also, results of 
  Hemke for a certain class of meromorphic functions 
  imply in the exponential case
  that, if the singular orbit of $\Ek$ satisfies
  a parabola condition, then the orbit of almost every point accumulates
  precicely on the postsingular set
  \cite[Corollary 6.1]{martinthesis}. As a final illustration of the
  use of our model in determining escape speed, we derive
  a combinatorial condition for this type of behavior.
  It is easy to see 
  that we need only consider the case of escaping endpoints, since
  any other 
  escaping point satisfies the parabola condition for all $K>0$
  (compare \cite[Proposition 4.5]{expescaping}).

 \begin{thm}[Endpoints that Satisfy a Parabola Condition]
  Let $\s\in\Sequ_0$ be a fast external address, and define
    \[ b := \limsup_{n\to\infty} F^{-(n-1)}(2\pi|s_n|). \]
  Set $t_n := T(\F^{n-1}(\s,\ts)) = t_{\sigma^{n-1}(\s)}$. 
  Then for every
   $K > 0$, the following two conditions are
  equivalent: 
  \begin{enumerate}
   \item There are $C_1 > 0$ and $n_0\in\N$
              such that $2\pi|s_n|< C_1\cdot t_n^K$ for $n\geq n_0$.
          \label{item:parabola}
   \item There are $C_2>0$ and $n_0$ such that
                  $2\pi|s_n|<C_2\bigl(F^{n-1}(b)\bigr)^K$ for $n\geq n_0$.
          \label{item:minimalpotentialcondition}
  \end{enumerate}
 \end{thm}
\begin{remark}
  $b$ is the \emph{minimal potential} as defined in
  \cite{expescaping}. 
\end{remark}
\begin{proof}
 Set $t_n^* := t_{\sigma^{n-1}(\s)}^*$. Note that
  $b$ is the limit of the 
  decreasing sequence $F^{-(n-1)}(t_n^*)$. 
  By Lemma
  \ref{lem:tsstar}, (\ref{item:parabola}) is equivalent to the
  same statement with $t_n$ replaced by $t_n^*$. 
  In particular, (\ref{item:minimalpotentialcondition}) implies
   (\ref{item:parabola}). 

  To prove the other direction, suppose that
   $2\pi|s_k|< C_3\cdot (t_n^*)^K$ for $n\geq n_0$. 
   Since $t_n^*\to\infty$,
   for sufficiently large $n$ we will get
   \[ 2\pi|s_n| < C_3\cdot (t_n^*)^K < \left(\frac{t_n^*}{K+1}+1\right)^{K+1}
                                                  -1,
     \]
   and similarly 
    \[ t_n^* < \left(\frac{t_n^*}{K+1}+1\right)^{K+1}-1. \]
   It follows that
   \[ t_{n}^* = F^{-1}\bigl(\max(t_{n+1}^*,2\pi|s_{n+1}|)\bigr) <
      F^{-1}\left(\left(\frac{t_n^*}{K+1}+1\right)^{K+1}-1\right) = 
             (K+1)F^{-1}\left(\frac{t_{n+1}^*}{K+1}\right). \]
   In other words, we have
    $t_n^*/(K+2) < F^{-1}(t_{n+1}^*/(K+2))$, which implies by induction that
   $t_n^* \leq (K+2)\cdot \bigl(F^{n-1}(b)\bigr)^K$. In particular,
    $2\pi|s_n| < C_3\cdot(K+2)^K\cdot \bigl(F^{n-1}(b)\bigr)^K$.
 \end{proof}


\section{Canonical Correspondence and Escaping Set Rigidity}
  \label{sec:uniqueness}
\label{sec:rigidity}

The bijection ${\DRconj}:X\to I(\Ek)$ 
 constructed in Section \ref{sec:classification} --- 
 while having certain continuity
 properties --- is, in general, quite far from being a conjugacy. 
 The question presents itself whether it is possible to construct
 a different map which is a conjugacy,
 or at least is continuous on a larger set.
 We will now show that this is not the
 case. 

 The underlying reason for this is the 
  rigidity presented by the $2\pi i$-periodic 
  structure of the dynamical plane: since the real and imaginary
  directions interact under an exponential maps, this rigid
  structure means that the escape speed
  of two points which correspond to each
  other under a conjugacy cannot dramatically differ.
  This idea is quite similar to that used by Douady and Goldberg
  \cite{douadygoldberg} showing that for $\kappa_1,\kappa_2\in
  (-1,\infty)$, the maps $E_{\kappa_1}$ and $E_{\kappa_2}$ are not
  conjugate on their Julia sets (and we obtain
  a generalization of their result in
  Theorem \ref{thm:noconjugacyescaping} below). 

 \begin{thm}[No Nontrivial Self-Conjugacies] \label{thm:noselfconjugacies}
  Let $Q>0$ and  
  suppose that
    $f: Y_Q\cap X \to X$
  is a continuous map with $f\circ \F = \F \circ f$ and 
    \begin{equation}
        f(\r,t) \in \{\r\}\times [0,\infty) \label{eqn:rayspreserved}
    \end{equation}
  for all $\r$ and $t$. Then $f$ is the identity.

  More precisely, let $(\s,t_0)\in Y_Q\cap X$. Suppose that a function
    \[ f: Y_Q \cap \bigcup_{j\geq 0} \F^{-j}(\F(\s,t_0)) \to X \]
   satisfies $f\circ \F = \F \circ f$ and (\ref{eqn:rayspreserved}).
    If $f(\s,t_0)\neq(\s,t_0)$, then $f$ is not continuous in $(\s,t_0)$.
 \end{thm}
 \begin{proof}
  The first statement follows immediately from the second.
  So let $x=(\s,t_0)$ and $f$ be as in the second part of the theorem,
  and suppose that $f(\s,t_0)=(\s,t_0')$ with $t_0\neq t_0'$.
  If $n$ is
  large enough, we can find $m(n)\in\N$ such that 
   \[ T(\F^n(x)) \leq \log(2\pi m(n)+t_0+1) < T(\F^n(x)) + 1. \]
  Let us define $y_n := \bigl(m(n)\, s_2\, s_3\, s_4 \dots,
  t_0\bigr)$. Now pull back the points $y_n$ along the orbit of
  $x$. More precisely, let $z_n$ be the uniquely defined point with address
  $s_1 s_2 \dots s_{n} m(n) s_2 s_3 \dots$ such that
  $\F^{n}(z_n) = y_n$.

  By the choice of $m(n)$,
   \begin{align*}
   T(\F^{n-1}(x)) &\leq T(\F^{n-1}(z_n))  < T(\F^{n-1}(x)) + 1,
      \quad\text{and thus} \\
   T(\F^j(x)) &\leq  T(\F^j(z_n))  < T(\F^j(x)) + F^{-(n-j-1)}(1)
   \end{align*}
  for every $j < n$. In particular, $z_n\in Y_{Q}$ and $z_n\to x$.
 
  On the other hand, consider the image points $f(z_n)$. Because
   $\F^n(f(z_n)) = f(y_n) = ( m(n) s_2 s_3 \dots , t_0')$, we see that
   \[ |T(f(z_n))-T(z_n)| < F^{-n}\bigl(|t_0'-t_0|\bigr) \to 0. \] 
  Thus
   $\displaystyle{\lim_{n\to\infty} f(z_n) = \lim_{n\to\infty} z_n =
         x \neq f(x) = f\left(\lim_{n\to\infty} z_n\right)}$,
  and $f$ is not continuous in $x$.
\end{proof}

\begin{cor}[No Other Conjugacies] \label{cor:nootherconjugacy}
 Suppose that $E_{\kappa_1}$ and $E_{\kappa_2}$ are exponential maps
 with nonescaping
 singular values which are conjugate on their
 sets of escaping points by a conjugacy $f$ that sends each dynamic
 ray of $E_{\kappa_1}$ to the corresponding dynamic ray of
 $E_{\kappa_2}$. Then $f={\DRconj}^{\kappa_2}\circ ({\DRconj}^{\kappa_1})^{-1}$.
\end{cor}
\begin{proof} Define $\Phi : I(E_{\kappa_1})\to I(E_{\kappa_1})$
   by $\Phi := f^{-1}\circ {\DRconj}^{\kappa_2}\circ
   ({\DRconj}^{\kappa_1})^{-1}$, and abbreviate ${\DRconj}:= {\DRconj}^{\kappa_1}$.
    If $Q$ is large enough, then $\Phi$ is continuous on
   ${\DRconj}^{\kappa_1}(Y_Q)$. Suppose that $z={\DRconj}(\s,t_0)\in I(E_{\kappa_1})$
   such that $\Phi(z)\neq z$. We may assume (by possibly exchanging
   $\kappa_1$ and $\kappa_2$ and passing to a forward image of $z$ if
   necessary) that $z\in Y_{Q}$ and $T({\DRconj}^{-1}(\Phi(z)))>t_0$. 
   Then the function 
    \[ f: \bigcup_{j\geq 0} \F^{-j}(\F(\s,t_0)) \cap Y_{Q} \to Y_Q,
      (\r,t) \mapsto {\DRconj}^{-1}(\Phi({\DRconj}(\r,t)))
    \]
   is continuous, contradicting Theorem
    \ref{thm:noselfconjugacies}.
 \end{proof}

The condition of every dynamic ray being sent to the
 corresponding ray in the limit dynamics should be satisfied by every
 ``reasonable'' conjugacy (up to a possible relabeling of the
 combinatorics).  
 A natural condition placed upon a topological conjugacy between two
 functions on an open set is that it preserves orientation. In our
 setting of conjugacies on the sets of escaping points, we replace
 this condition by a notion of ``order-preserving'' conjugacies.

 The collection of dynamic rays is endowed with a natural vertical
 order: of any two dynamic rays,
 one is \emph{above} the other. More precisely, 
 define $\Hplane_R := \{z\in\C: \re z > R\}$. 
 If $\gs$ is a dynamic ray and $R$
 is large enough, then  the set $\Hplane_R \setminus
 \gs\bigl([1,\infty)\bigr)$ has exactly two
 unbounded components,
 one above and one below $\gs$, and any other dynamic ray
 must tend to
 $\infty$ within one of these. It follows immediately from the
 construction of dynamic rays that this order coincides precisely with
 the lexicographic order of their external addresses.

 We now call a continuous function
 from some subset of $X$ to $X$ \emph{order-preserving} if it induces
 an order-preserving map
 on external addresses. Similarly,
 if $\kappa_1,\kappa_2\in\C$ and $f:I(E_{\kappa_1})\to
 I(E_{\kappa_2})$ is continuous, then we call $f$ order-preserving if
 it preserves the vertical (and thus the lexicographic) order of
 dynamic rays. Any orientation-preserving
 conjugacy of two exponential maps induces an order-preserving map
 on their sets of escaping points. 

 There are only a few order-preserving self-conjugacies of the
  shift map on external addresses. The simple proof of the following fact
  is left to the reader

\begin{lem}[Self-Conjugacies of the Shift] \label{lem:shiftconjugacies}
 Let $f$ be an order-preserving homeomorphism of the space of external addresses such
  that $f\circ\sigma = \sigma\circ f$. Then there exists $j\in\Z$ such that
  \[ f(s_1 s_2 s_3 \dots ) = (s_1 + j) (s_2 + j) (s_3 + j) \dots \]
 for all external addresses $\s=s_1 s_2 s_3 \dots$ . \qedd
\end{lem}

Using Corollary \ref{cor:nootherconjugacy}, we can now describe all
 possible order-preserving conjugacies between the escaping dynamics
 of two exponential maps.

 \begin{cor}[Order-Preserving Conjugacies]
  Suppose that $f$ is an order-preserving conjugacy between two maps
  $E_{\kappa_1}$ and $E_{\kappa_2}$ (with nonescaping singular orbit)
  on their escaping sets. 
  Then there exists a parameter
  $\kappa_2' = \kappa_2 + 2\pi i k$ (with $k\in\Z$) such that
  \[ {\DRconj}^{\kappa_2'}\circ ({\DRconj}^{\kappa_1})^{-1}: 
        I(E_{\kappa_1})\to
  I(E_{\kappa_2'}) \] 
  is a conjugacy.
 \end{cor}
\begin{proof}
  The map $f$ induces an order-preserving self-conjugacy of the shift. By Lemma
  \ref{lem:shiftconjugacies} this map consists of shifting all labels
  by some number $k$. Let $\kappa_2' := \kappa_2 - 2\pi i k$.
  The maps $E_{\kappa_2}$ and $E_{\kappa_2-2\pi ik}$ are conjugate by
  the map $z\mapsto z - 2\pi i k$, and the induced self-conjugacy of the
  shift consists of shifting all labels by $-k$. Thus the map
  $f':z\mapsto f(z)-2\pi i k$ is a conjugacy between $E_{\kappa_1}$
  and $E_{\kappa_2}$ which preserves dynamic rays. The claim follows
  by Corollary \ref{cor:nootherconjugacy}. 
\end{proof}

The next theorem is a generalization of the previously mentioned
 result of Douady and
 Goldberg \cite{douadygoldberg}.

\begin{thm}[No Conjugacy for Escaping Parameters] \label{thm:noconjugacyescaping}
 Let $\s\in\Sequ_0$ and let
 $(\s,t_1),(\s,t_2)\in X$ with $t_1\neq t_2$.
 Suppose that $\kappa_1$ and $\kappa_2$ are parameters such that 
 ${\DRconj}^{\kappa_1}(\s,t_1) = \kappa_1$ and ${\DRconj}^{\kappa_2}(\s,t_2) =
 \kappa_2$. Then $E_{\kappa_1}$ and $E_{\kappa_2}$ are not conjugate
 on $\C$.
\end{thm}
\begin{proof}
  By contradiction, let $f:\C\to\C$ be a conjugacy between
 $E_{\kappa_1}$ and $E_{\kappa_2}$. For some $Q>0$, the map
  \[ \alpha: Y_Q\to X; (\s,t)\mapsto 
      \left({\DRconj}^{\kappa_2}\right)^{-1}(f({\DRconj}^{\kappa_1}(\s,t))) \]
 is defined. $\alpha$ is either order-preserving or order-reversing,
 depending on whether $f$ is orientation-preserving or
  -reversing. Also $\alpha\neq\id$ because 
  $\alpha(\F^n(\s,t_1)) = \F^n(\s,t_2)\neq \F^n(\s,t_1)$. 

 Suppose first that $\alpha$ is order-preserving. Since
  $f$ must map $g_{\sigma^n(\s)}^{\kappa_1}$ to 
  $g_{\sigma^n(\s)}^{\kappa_2}$, it follows by Lemma
  \ref{lem:shiftconjugacies} that $\alpha$ preserves 
  external addresses (i.e. satisfies
  (\ref{eqn:rayspreserved})). 
  As in the proof of Corollary \ref{cor:nootherconjugacy},
  this contradicts
  Theorem \ref{thm:noselfconjugacies}.

 Now suppose that $\alpha$ is order-reversing. Since $f$ maps
  postsingular rays of $E_{\kappa_1}$ to those of $E_{\kappa_2}$,
  this is possible only when $\s$ is periodic of period $1$.
  We may assume without loss of generality that 
  $\s=\per{0}$, in which case we can replace $\alpha$ in the above
  argument
  by the order-preserving map $\wt{\alpha}(s_1 s_2 \dots) :=
       \alpha((-s_1) (-s_2) \dots )$.
\end{proof}

We are now in a position to extend Proposition \ref{prop:noconjugacy} to a
 larger class of examples.

\begin{thm}[Rigidity for Parameters with Different Combinatorics]
  \label{thm:rigidity}
 Suppose that $\s\in\Sequ_0$ and $\wt{\kappa}$ 
 is a nonescaping parameter for which the singular value
 is contained in the limit set of
 $g^{\wt{\kappa}}_{\s}$. 
 Suppose furthermore that $\kappa$ is another
 nonescaping parameter such
 that the limit set of $g^{\kappa}_{\s}$ does not contain the singular value
 and is bounded. Then 
 ${\DRconj}^{\wt{\kappa}}\circ ({\DRconj}^{\kappa})^{-1}$ is not
 continuous.
\end{thm}
\begin{proof}
 Suppose by contradiction that 
  ${\DRconj}:= {\DRconj}^{\wt{\kappa}}\circ ({\DRconj}^{\kappa})^{-1}$ is
 continuous. As in the 
 proof of Proposition \ref{prop:noconjugacy}, pick an arbitrary
 point $w\in I(E_{\kappa})$ and a neighborhood $U$ of $w$ whose
 image under ${\DRconj}$ is bounded. 

 Let $A$ denote the accumulation set of $g^{\kappa}_{\s}$. We shall
  show that there is an iterated preimage of $A$ which is contained in 
  $U$. The conclusion then follows in the same way as in Proposition
  \ref{prop:noconjugacy}. 

First note that every component of the
 preimage of $A$ is compact. Indeed, otherwise there is no continuous
 branch of $\Ek^{-1}$ on $A$, which means that $A$ separates $\kappa$
 from $\infty$. However, this is impossible: if $\kappa\in J(E_{\kappa})$,
 then there must be escaping points close to $\kappa$, which are
 connected to $\infty$ by a dynamic ray, and if $\kappa\in F(E_{\kappa})$,
 then it is easy to see that there is a curve in the Fatou set which
 connects $\kappa$ to
 $\infty$ (see \cite{expattracting} or
 Section \ref{sec:topologyattracting}). All
 preimages of $A$ are translates of each other; let $K$ denote the
 diameter of any of these preimages.
 Choose $n$ sufficiently large and let $V\subset U$ be
 a small neighborhood of $w$ which is mapped biholomorphically to
 $\D_{2\pi+1}(\Ek^n(w))$ by $\Ek^n$. (The existence of such a $V$ is
 easily shown using a pullback argument.) Chose among the preimages
 of $A$ one, call it $A_0$, which satisfies 
   \[ |\Ek^{n+1}(w)-\kappa| \leq |z-\kappa| \leq |\Ek^{n+1}(w)-\kappa| + 2\pi + K \]
  and for all $z\in A_0$. If $n$ was chosen large enough,
  then $|\Ek^{n+1}(w)-\kappa| > 2\pi + K$ and thus
  \[ \log|z-\kappa| - \log |\Ek^{n+1}(w)-\kappa| \leq 
     \log\left(1+\frac{2\pi+K}{|\Ek^{n+1}(w)-\kappa|}\right) \leq
      \frac{2\pi + K }{|\Ek^{n+1}(w)-\kappa|} < 1.\]
 Thus if we take the pullback $A_1$ of $A_0$ by the same branch of
 $\Ek^{-1}$ that carries $\Ek^{n+1}(w)$ to $\Ek^n(w)$, then $A_1\subset
 \D_{2\pi+1}(\Ek^n(w))$. We can then further pull back $A_1$ to
 $V\subset U$, which concludes the proof.
\end{proof}

The
preceding result, together with theorems on the combinatorial rigidity
of escaping and Misiurewicz parameters
\cite{markusdierk,markuslassedierk,expmisiurewicz}, 
 easily implies the following statement on the rigidity of the
 escaping dynamics of such parameters.

\begin{cor}[Rigidity for Escaping and Misiurewicz Parameters]
 Suppose that $\kappa_1\neq \kappa_2$ are attracting, parabolic,
 Misiurewicz or escaping parameters, at
 least one of which is not attracting or parabolic. Suppose that
  $\im\kappa_1,\im\kappa_2\in (-\pi,\pi]$. Then
 $E_{\kappa_1}$ and $E_{\kappa_2}$ are not conjugate on their sets of
 escaping points by an order-preserving conjugacy. 
\end{cor}
\begin{proof} Clearly an escaping parameter cannot be conjugate to a
 nonescaping parameter. So let us first suppose that $\kappa_1$ and
 $\kappa_2$ are escaping parameters, and that the singular values
 lie on the rays at external addresses $\s^1$ and $\s^2$. Then
 both addresses have first entry $0$ by 
 \cite[Corollary 1]{markuslassedierk}. 
 Since the conjugacy must map the
 singular value of $E_{\kappa_1}$ to that of $E_{\kappa_2}$, it follows
 by Lemma \ref{lem:shiftconjugacies} that $\s^1=\s^2$. 
 As in the proof of Theorem \ref{thm:noconjugacyescaping},
 their potentials are equal as well. However, this contradicts 
 the fact that for every $x\in X$ there exists only one parameter
 $\kappa$ with $\DRconj^{\kappa}(x)=\kappa$ 
 \cite[Theorem 2]{markuslassedierk}. 

 Now suppose that both $\kappa_1$ and 
 $\kappa_2$ are Misiurewicz. Assume that the preperiod of
 $\kappa_1$ is smaller or equal to that of $\kappa_2$. By 
 \cite[Theorem 4.3]{expper}, 
 there exists a preperiodic address $\s$ such that $\gs^{\kappa_1}$
 lands at $\kappa_1$.  By \cite[Theorem 3.2]{expper}, all periodic rays
 of $E^{\kappa_2}$ land at periodic points.
 Because the preperiod of $\kappa_2$ is greater or equal to that of
 $\kappa_1$, this implies that $\gs^{\kappa_2}$ lands at a preperiodic
 point. 
 By the results of \cite{expmisiurewicz}, Misiurewicz parameters with
 given combinatorics are unique, so this landing
 point is $\neq \kappa_2$ . Thus we can apply
 Theorem \ref{thm:rigidity}. The same argument works (without
 reference to \cite{expmisiurewicz}) if $\kappa_2$ is parabolic
 or attracting. 
\end{proof}

It seems reasonable to conjecture that the escaping dynamics of exponential
maps whose singular value lies in the Julia set is always rigid.
This conjecture would imply density of hyperbolicity: a non-hyperbolic
stable parameter would be (quasiconformally) conjugate to all nearby
parameters, and in particular the maps would be conjugate on their
sets of escaping points.

  \begin{conj}[Escaping Set Rigidity]
   Suppose that $\kappa_1$ is a parameter with $\kappa_1\in
   J(E_{\kappa_1})$, and let $\kappa_2\notin \{\kappa_1+2\pi i k\}$. Then there
   exists no order-preserving conjugacy
    \[f:I(E_{\kappa_1})\to I(E_{\kappa_2})\]
   between $E_{\kappa_1}$ and $E_{\kappa_2}$.
  \end{conj}

 \section[Topology of the Julia Set]{Topology of the Julia Set for
    Attracting and Parabolic Parameters}
    \label{sec:topologyattracting}

  We will now completely describe the Julia sets of
  attracting
  and parabolic exponential maps (and the dynamics thereon)
  as a quotient
  of our model $\Xb$. In particular, any two attracting exponential maps are
  conjugate on their sets of escaping points. We will give the
  complete construction for attracting parameters, and remark on the
  parabolic case later.

  So let $\kappa\in\C$ such that $\Ek$ has an attracting cycle
  $a_0\mapsto a_1 \mapsto \dots \mapsto a_n=a_0$ and corresponding
  Fatou components $A_0\mapsto A_1 \mapsto \dots \mapsto A_n = A_0$.
  This cycle of immediate attracting basins must contain the singular
  value \cite[Theorem 7]{waltermero};
  let us choose our labelling
  in such a way that $\kappa\in A_1$.
  By the Koenigs linearization theorem 
  \cite[Theorem 8.2]{jackdynamics}, we can find 
  open Jordan neighborhoods
  $V_j$ of $a_j$ such that $\kappa\in U_1$, $\Ek(V_j)\compin
  V_{j+1}$%
  \footnote{The notation $U\compin V$ means, as usual, that $\cl{U}$ is a
  compact set contained in $V$}
  and $\Ek(U_n)\compin U_1$. For $j=0,\dots,n-1$,
  let $U_j$ denote the component of
  $\Ek^{-1}(V_{j+1})$ containing
  $a_j$. Since
  $\kappa\in U_1$, the component $U_0$ contains a left halfplane;
  the other $U_j$ are bounded Jordan domains. 
  We consider the set
    \[ W := \C\setminus\left(\bigcup_{i=0}^{n-1} \cl{U_i}\right). \]
  Then $\cl{\Ek^{-1}(W)}\subset W$, and $\Ek:\Ek^{-1}(W)\to W$ is a
  covering map. 

 From now on let us suppose that $n\geq 2$. The minor modifications
  necessary in the
  case $n=1$ are straightforward and are left to the reader.

 We can connect $\kappa$ and $E^n(\kappa)$ by a curve in $U_1$
  (e.g.~by a straight line in linearizing coordinates). Pulling this
  curve back under $\Ek^n$, we obtain a curve 
  $\gamma\subset A_1$ which connects $\kappa$ to $\infty$.
  Define
  $V:=\C\setminus\gamma$. 
  The set
  $\Ek^{-1}(V)$ then consists of countably many strips bounded by two
  preimages of $\gamma$.
  Let us label these strips as $R_k$ in such a
  way that $t+(2k+1)\pi i \in R_k$ for large enough $t$, and let
  \[ \wt{L}_k:V\to R_k \]
  denote the corresponding branch of $\Ek^{-1}$.
  Observe that these 
   differ from the branches $L_k$ considered in Section
   \ref{sec:classification}. Note that
   $\wt{L}_k$ is well-defined everywhere on the Julia set.

 If $z\in J(\Ek)$, we can associate to $z$ an \emph{itinerary}
   $\itin(z) := \u_1 \u_2 \u_3\dots$ such that
   \[ \Ek^{n-1}(z) \in R_{\u_n} \]
  for all $n\geq 1$. If $\s\in\Sequ_0$, then
   all points in $\gs$ clearly have the same itinerary, which is also
   denoted by $\itin(\s)$. (This itinerary can be defined in a purely
   combinatorial way by associating an ``intermediate external
   address'' to the curve $\gamma$; see \cite{expattracting,expper}
   or, for a general approach, \cite[Section 3.7]{thesis}.)

 Choose some $A>0$  and $B<0$ such that the map
  ${\DRconj}={\DRconj}^{\kappa}:X\to I(\Ek)$ satisfies
  $|\re {\DRconj}(\s,t) - t | < 2$ on $Y_A$ (as in Theorem
  \ref{thm:conjugacyconvergence}) and such that
   \[
      \Hplane := \{z\in\C: \re z > A-2 \} \subset W \subset
      \{z\in\C: \re z > B \}. \]
 
 Now let us define functions $H_k:\Xb\to J(\Ek)$ 
  by
  \[ H_0(\s,t) := {\DRconj}(\s,t+A)\quad \text{ and }\quad
    H_{k+1}(\s,t) := \wt{L}_{\u_1}(H_k(\F(\s,t))), 
  \]
  where $\u_1$ is the first entry of
  $\itin(\s)$. Note that
  $H_k(\s,t)$ always lies on the dynamic ray $g_{\s}$.

 \begin{thm}[Conjugacy for Attracting Parameters] \label{thm:topologyattractingmain}
  In the hyperbolic metric of $W$, the functions $H_k$ converge
   uniformly to a continuous, surjective function
   $H: \Xb \to J(\Ek) $
  with $H\circ \F = \Ek \circ H$.
  Furthermore, $H|_X: X\to I(\Ek)$ is a conjugacy.
 \end{thm}
\begin{remark}
  By Corollary \ref{cor:nootherconjugacy}, 
   $H_{|X}$ must be equal to ${\DRconj}$. In particular, 
   Theorem 9.1 implies Theorem \ref{thm:attractingconjugacy}.
\end{remark}
\begin{proof}
 Let us denote the 
  hyperbolic metric of any hyperbolic domain $U\subset\C$ by $ds = \rho_U |dz|$.
  Since $\Ek:\Ek^{-1}(W)\to W$ is a covering map, it expands the
  hyperbolic metric of $W$. In fact, there exists $K>1$ such that
  $\|D\hspace{-0.5mm}\Ek(z)\|_{\operatorname{hyp}} \geq
  K$ for all $z\in
  \Ek^{-1}(W)$.

 To prove this, set $W' := \Ek^{-1}(W)$, and let
  \[ \wt{W} := \bigl\{z\in W: z+2\pi i k \in W \text{ for all
  $k\in\Z$}\bigr\}.\]
  (Thus $\wt{W}$ is obtained from $W$ by removing all translates of
   the sets $\cl{U_i}$.) Then $W'\subset \wt{W}\subset W$.
  Because $\rho_{\wt{W}} \geq \rho_{W}$ by monotonicity of the
  hyperbolic metric, it is sufficient to show
  that, for every $z\in W'$, 
  \begin{equation}
    q(z) := \frac{\rho_{W'}(z)}{\rho_{\wt{W}}(z)} \geq K >1.
  \label{eqn:hyperbolicmetricquotient} \end{equation} 
  Recall that $U_0$ contains a left halfplane, so that, for
  some $R_0>0$,
  the set $\C\setminus W'$ contains the curves 
  $\{R+(2k+1)\pi i:R\geq R_0\}$. By 
  standard estimates on the hyperbolic metric, it
  follows that $\rho_{W'}(z)$ is bounded from below as $\re z \to +\infty$. 
  On the other hand, 
  $\wt{W}$ contains the right halfplane $\Hplane$, so
  $\rho_{\wt{W}}(z)\to 0$, and thus
  $q(z)\to\infty$, as $\re z\to +\infty$.
  Since $\cl{W'}\subset \wt{W}$ and 
  $\wt{W}$ is bounded to the left, it follows 
  that 
   \[ K := \inf_{\genfrac{}{}{0pt}{}{z\in W':}{|\im z|\leq \pi}} q(z) > 1. \]
  The expression $q(z)$ is $2\pi i$-periodic, so
  (\ref{eqn:hyperbolicmetricquotient}) follows.

 For an arbitrary $(\s,t)\in X$, consider the two points
  $z_1 := \Ek(H_0(\s,t)) = {\DRconj}(\F(\s,t+A))$ and
  $z_2 := \Ek(H_1(\s,t)) = H_0(\F(\s,t))$. Both points have real
  parts greater than $A-2$ and thus can be connected by a straight line
  $g_0$ in $W$. Note that $g_0$ is homotopic (in $W$) to the piece of
  the ray $g_{\s}$ between $z_1$ and $z_2$, as this piece is also
  contained in the halfplane $\Hplane$. Thus we can pull back
  $g_0$ and obtain a curve $g_1$ between $H_0(\s,t)$ and
  $H_1(\s,t)$. We claim that the (euclidean) length of $g_1$ is
  uniformly bounded (independent of $\s$ and $t$).
 
 To prove this claim, recall that 
   \[ \re(z_1) \leq T(\F(\s,t+A)) + 2 = \exp(t+A)-2\pi|s_2| + 1 \]
  and
   \[ \re(z_2) \geq T(\F(\s,t)) - 2 = \exp(t)-2\pi|s_2| - 3. \]
 It follows that the euclidean length of $g_0$ satisfies
  \[ \ell(g_0) \leq \exp(t+A)-\exp(t)+4+2\pi = O(\exp(t)). \]
 Because all points of $g_0$ have absolute value at least
  $|z_2|-2\pi \geq \frac{\exp(t)}{\sqrt{2}} - 2 - 3\pi$, we see that
  \begin{align*}
    \ell(g_1) \leq \frac{1}{|z_2|-2\pi}\cdot \ell(g_0) = O(1).
  \end{align*}
 
 Since $\rho_W(z)\to 0$ as $\re z\to\infty$, the function
  $\rho_W$ is uniformly bounded on $W'$. Thus the hyperbolic length of
  $g_1$ in the hyperbolic metric of $W$ is also bounded by some
  constant $C$. Now, taking pullbacks inductively, we see that the
  hyperbolic distance between $H_k(\s,t)$ and $H_{k+1}(\s,t)$ is
  bounded by $\frac{C}{K^k}$. Thus the $H_k$ converge uniformly. The
  functional equation $H\circ \F = \Ek \circ H$ is satisfied by construction.
 
 To show surjectivity of $H$, it is sufficient to see that $H(\Xb)$ is
 dense in $J$ (note that, because the hyperbolic distance between $H(x)$ and
 $H_0(x)$ is uniformly bounded, $H$ is
 again continuous as a map $\Xb\cup\{\infty\}\to
 J\cup\{\infty\}$). However, density of the image is trivial because
 $\Ek^{-1}(H(\Xb))\subset H(\Xb)$, and backward orbits of any point
 (except $\kappa$) accumulate on the entire Julia set. 
 Injectivity of $H$ on $X$ follows by the same argument as before.
\end{proof}

The following immediate corollary was previously proved in
\cite{accessible} (with somewhat different notation).
\begin{prop}[Dynamic rays landing at a common point]
   \label{prop:raysitinerary}
 Let $\kappa$ be an attracting parameter. Then every
  non-escaping point in $J(\Ek)$ is the landing point of at least one
  dynamic ray. Two dynamic rays land at the same point if and only if
  they have identical itineraries.
\end{prop}
\begin{proof} The only thing we still need to check is that, whenever
 $\itin(\s) = \itin(\s')$, then $H(\s,\ts) = H(\s',t_{\s'})$. However,
 the strips $R_k$ have height $2\pi$, so
 the $n$-th entries of $\s$ and $\s'$, for any $n$, differ by at
 most $1$. It follows easily (for example by Lemma
 \ref{lem:tsstar}) that
 $|t_{\sigma^{n}(\s)}-t_{\sigma^n(\s')}|$ is bounded independently of $n$. Thus
 the distance between
 the points $H_0(\sigma^n(\s),t_{\sigma^n(\s)})$ and
 $H_0(\sigma^n(\s'),t_{\sigma^n(\s')})$ is uniformly bounded
 as well. The claim now follows by the
 contraction argument from the previous proof.
\end{proof}

Note that the proof of \ref{prop:raysitinerary} for periodic
 addresses is much easier, compare \cite[Proposition 4.5]{expper}. 

 Another direct consequence of Theorem
  \ref{thm:topologyattractingmain} is the following result, which
  describes an abstract model of the Julia set of an attracting
  exponential map, in analogy to the ``Pinched Disk Model'' for
  polynomials \cite{pincheddisk}.

\begin{cor}[``Pinched Cantor Bouquet'']
 Let $\kappa$ be an attracting parameter.
 Form the quotient $\wt{X}$ of $X$ by identifying all points
  $(\s,\ts)$ and $(\s',t_{\s'})$ for which 
   $\itin(\s) = \itin(\s')$.
 Then $\F$ projects to a map $\F:\wt{X}\to\wt{X}$ which is conjugate
 to ${\Ek}:J(\Ek)\to J(\Ek)$.
 \ep
\end{cor}

\begin{remark}
 All the preceding theorems remain true 
 for parabolic parameters. The issue is to find a
 replacement for the strict hyperbolic contraction used in the proof of
 Theorem \ref{thm:topologyattractingmain}.
 This issue is the same which appears in the proof of local
 connectivity for quadratic polynomials with a parabolic orbit --- see
 \cite[Expos\'e 10]{orsay} or 
 \cite{geometricallyfinitelocalconnectivity} --- and can be dealt with in
 a similar manner.
 Those arguments, however, are somewhat technical and hardly very
 enlightening; it seems to us that there is
 little to be gained by their detailed adoption to the exponential case.
 Furthermore, in recent work by Haissinsky
 \cite{haissinsky}, parabolic rational maps were constructed from
 hyperbolic maps by using Guy David's
 transquasiconformal surgery. In particular, the resulting parabolic
 map is topologically conjugate to the hyperbolic function it
 originated from. Such methods should generalize to the
 space of exponential maps and thus yield a natural proof of the
 conjugacy of parabolic exponential maps to attracting exponential
 maps with the same combinatorics (i.e., with the same \emph{intermediate
 external address} \cite{expattracting}). 
 In view of these facts, we have decided against a
 presentation of rigorous proofs of the above theorems in the
 parabolic case. 
\end{remark}

\section{Invalidity of Renormalization} 
\label{sec:norenormalization}

Suppose that $\Ek$ is any attracting exponential map of
period $n>1$ and let $\mu$ be the multiplier of its attracting orbit.
As in the previous section, label the cycle of immediate basins,
  $A_0\mapsto A_1\mapsto\dots\mapsto A_n=A_0$ in such a way that 
 $A_0$ contains a left
half plane. $\Ek^n|_{A_0}$ is conformally conjugate to
$E_{\kappa_0}|_{F(E_{\kappa_0})}$ where $\kappa_0$ is such that
$E_{\kappa_0}$ has an attracting fixed point with multiplier
$\mu$ (in fact, $\kappa_0=\log\mu - \mu$). This can be proved
either by constructing the conjugacy directly using the linearizing
coordinates of $\Ek$ and $E_{\kappa_0}$, or by conjugating these maps
to a normal form as in 
 \cite[Section III.4]{habil} or \cite{devgoldberg}. 
Let
 \[ \Psi: A_0 \to F(E_{\kappa_0}) \]
be this conjugacy, and note that
  $\Psi(z+2\pi i)=\Psi(z)+2\pi i$. 
We will now prove Theorem \ref{thm:renorm}, i.e.\ that this map does not
 extend continuously to $\cl{A_0}$.
 The reason, as in the argument from Theorem
 \ref{thm:noselfconjugacies}, is that the $2\pi i$-periodic structure
 of the dynamical plane must be preserved under a conjugacy, which
 makes it impossible to conjugate $E_{\kappa_0}$ to the much faster
 growing function $\Ek^n$. Combinatorially speaking, this means that
 there are dynamic rays on $\partial A_0$ which would be mapped to
 points with an external address which is not exponentially bounded,
 which is clearly impossible (compare the combinatorial tuning formula
 in \cite{expcombinatorics}). However, our proof does not use these
 combinatorial notions and is, in fact, completely elementary.

\begin{thm}[No Topological Renormalization] \label{thm:norenormalization}
 Let $\Ek$ and $\Psi$ be as above. Then $\Psi$ does not have a
 continuous extension to $\partial A_0$.
\end{thm}
\begin{proof} Assume, by contradiction, that $\Psi$ does extend continuously
 to $\partial A_0$. The idea of the proof is the following: orbits
 of points in $\cl{A_0}$ under $\Ek^n$ (with bounded imaginary parts)
 grow essentially as iterates of $F^n$. So for any large enough $K$
 and $k$,  we can find a point $z_0$ with
 real part around $K$ whose imaginary part under
 $\Ek^{kn}$ is about $F^{kn}(K)$. The orbit of $\Psi(z_0)$, on the
 other hand, can grow only like $F$ under iteration of $E_{\kappa_0}$,
 which leads to a contradiction because
 $\Psi(\Ek^{kn}(z_0))=E_{\kappa_0}^k(\Psi(z_0))$ must also have
 imaginary part roughly $F^{kn}(K)$. In the following, we fix the
 details of the proof.

 Cut the plane into the 
 strips $R_k$ from the previous section; recall that these strips
 have bounded imaginary
 parts and the strip boundaries lie in $A_0$. 
 We may assume that the $R_k$ are numbered so that $R_0$ contains $0$. Let 
   \[ R := R_0 \cup \bigcup_{1\leq j \leq n-1} R_{\u_j}, \]
  where $R_{\u_j}$ is the strip containing $A_j$. 
 Then the orbit of any point $z\in\cl{A_0}$ lies in
 $\cl{A_0}\cup R$. Define $A := \max_{z\in \cl{R}} |\im(z)|$ and
 $B:=\max_{t\in[-3\pi,-\pi]}|\im \Psi(ti)|$. Note that both quantities
 are at least $\pi$.

 Choose $K>B$ large enough such that, 
  whenever $|\Ek(z)| \geq K$, then
   \begin{align}  
    |\Ek(z)|-A-1&\geq F(\re(z)-1) \label{eqn:norenormalization1}
                \quad\text{and} \\
     |\Ek(z)|+A+1 &\leq F(\re(z)+1). \label{eqn:norenormalization2}
   \end{align}

 Let 
   \[ M := \max \left\{ |\Psi(z)|: z\in \cl{A_0\cap R_0} \text{ and } \re(z)
      \in [K-1,K+1] \right\}; \]
 by enlarging $M$, if necessary, we can also assume that 
 $\exp(t)+|\kappa_0|+1 \leq F(t+1)$ for all $t\geq M$. Note that this
   implies $|E_{\kappa_0}^j(z)| < F^j(M+1)$ for all $j$ and
   all $z$ with $\re(z) \leq M$. 
 Finally, choose $k$ so large that $F^{k}(K-B) > M +1 $. 

Pick any point $z_1\in
 \cl{A_0}$ with $\re(z_1) = 0$ and
  \[ \im(z_1) \in \left[F^{kn}(K)-\pi, F^{kn}(K)
   + \pi  \right]. \] 
 By repeatedly pulling back the point $z_1$ under $\Ek^{-n}$, we obtain a
 point $z_0\in
 \cl{R_0\cap A_0}$ with $\Ek^{kn}(z_0)= z_1$ and $\Ek^{jn}(z_0)\in \cl{R_0}$ 
  for $j < k$. By (\ref{eqn:norenormalization1}), we see that
  \[ \re\bigl(\Ek^{kn-1}(z_0)\bigr) -1 \leq 
     F^{-1}( |z_1| - A - 1 ) \leq F^{kn-1}(K) \]
 and, similarly, by (\ref{eqn:norenormalization2}),
  $\re\bigl(\Ek^{kn-1}(z_0)\bigr)+1 \geq
     F^{kn-1}(K)$.
 In particular, 
  \[ \re\bigl(\Ek^{kn-1}(z_0)\bigr) \geq K. \]
  Repeating this argument inductively, it follows that
  $\re( z_0) \in [K-1,K+1]$.
 
 Because $\Psi(z+2\pi i)=\Psi(z)+2\pi i$, we can estimate that
  $\im(\Psi(z_1)) \geq F^{kn}(K)-B$.
 On the other hand, $\re(\Psi(z_0)) \leq M$, and thus 
  \[ 
    |\Psi(z_1)| = |E_{\kappa_0}^k(\Psi(z_0))| < 
      F^k(M+1) <
      F^{2k}(K-B)
     \leq F^{kn}(K) - B. 
  \]
 This is a contradiction. 
\end{proof}

\section{Parameter Space} \label{sec:parameterspace}

 In \cite{markusdierk}, it was shown that the parameters
  for which the singular value
  lies on a dynamic ray (but is not the endpoint of a ray)
  are organized in \emph{parameter rays}. 

 \begin{prop}[Classification of Parameters on Rays
              \protect{\cite{markusdierk}}]
 \label{prop:parameterclassification}
  For every $\s\in\Sequ_0$ and
  every $t>\ts$, there exists a unique parameter $\kappa =
  {\PRconj}(\s,t)$ 
  such that ${\DRconj}^{\kappa}(\s,t)=\kappa$. For any fixed $\s$, 
  the value
  ${\PRconj}(\s,t)$ depends
  continuously on $t$ and satisfies
   $\im \PRconj(\s,t) \to 2\pi i s_1$ as
   $t\to\infty$. \qedd
 \end{prop}

The curves
   $\PR:(0,\infty)\to\C; t\mapsto {\PRconj}(\s,t+\ts)$,
  are called \emph{parameter rays}. 
  As for dynamic rays,
  we say that $\PR$ \emph{lands} at a parameter $\kappa\in\Ch$ if
  $\lim_{t\to 0} \PR(t)=\kappa$.
  If $\kappa\neq\infty$, we set
   we set $\PR(0) := \PRconj(\s,\ts) := \kappa$ in this case.
   (Similarly, in the following we will write
   $\DRconj^{\kappa}(\s,\ts)=z$ if the dynamic ray $\gs^{\kappa}$ lands
   at $z\in\C$.)

 In this section, we investigate continuity properties
  of the map $\PRconj$, using the results of Section
  \ref{sec:classification}. The means to transfer this
  dynamical information into the parameter plane is provided
  by the following result. (Recall the definition of
  $Q(K)$ and $Y_Q$ from (\ref{eqn:Q})). 

 \begin{lem}[Continuity of $\PRconj$] \label{lem:localcontinuity}
  Let $\kappa_0\in\C$. Suppose that
   there are $n\geq 0$, $Q_1 > Q(|\kappa_0|)$ and
   $x_0\in \F^{-n}(Y_{Q_1})$ such that $\DRconj^{\kappa_0}(x_0)=\kappa_0$.
   Then $\kappa_0 = \PRconj(x_0)$; furthermore
   there is a neighborhood
   $V$ of $x_0$ in $\F^{-n}(Y_{Q_1})$ such that the map
   $\PRconj:V\to \C$ is defined and a homeomorphism onto its image. 
 \end{lem}
 \begin{proof} Pick a neighborhood $U$ of $\kappa_0$ with
   $Q(|\kappa|)<Q_1$ for all $\kappa\in U$. If $V_1$ is a small
   neighborhood of $x_0$ in $\F^{-n}(Y_{Q_1})$ and
   $U$ was chosen small enough, then 
   $\DRconj^{\kappa}(x)$ is defined for all $\kappa\in U$ and
   $x\in V_1$, and jointly continuous in $\kappa$ and $x$. 

  By Hurwitz's theorem, there is a compact
   neighborhood $V\subset V_1$ of $x_0$
   and a function $V\to U; x\mapsto \kappa(x)$ such that
   $\DRconj^{\kappa(x)}(x)=\kappa(x)$. Furthermore, this map can be chosen
   in  such a way that
   $\kappa(x)\to \kappa_0$ as $x\to x_0$. 

  Suppose that $x_0=(\s,t_0)$. Then there is
   $t_1$ such that $(\s,t)\in V$ for
   all $t \in (t_0,t_1)$. By 
   Proposition \ref{prop:parameterclassification}, we have
   $\kappa(\s,t) = \PRconj(\s,t)$, and hence
   $\kappa_0 = \lim_{t\to t_0} \kappa(\s,t) = 
                 \lim_{t\to t_0} {\PRconj}(\s,t)$.
   Thus $\kappa_0 = \PRconj(\s,t_0)$, as required.

  For every $x\in V$, the parameter $\kappa(x)$ also satisfies the
   hypotheses of the theorem. Thus, by what
   we have just shown, $\kappa(x)=\PRconj(x)$ for
   all $x\in V$. So $\PRconj$ is defined on $V$ and continuous
   in $x_0$. Continuity in any other point of $V$ follows by 
   replacing $\kappa_0$ by $\PRconj(x)$.
   Since $\PRconj$ is clearly injective, and $V$ was chosen
   to be compact, $\PRconj|_{V}$ is a homeomorphism onto its image.
  \end{proof}

 \begin{cor}[Continuity away from Endpoints]
  For every $\eps>0$, the map ${\PRconj}$ is a homeomorphism when restricted to
   the set
  \[ Z_{\eps} := \{(\s,t):t\geq \ts+\eps \}. \]
 \end{cor}
 \begin{proof}
 Let $x_0\in Z_{\eps}$ and $\kappa_0 := {\PRconj}(x_0)$. 
  Choose $n$ sufficiently large such that
  $F^n(\eps) \geq Q(|\kappa|)+1$. Then
  $\F^n(Z_{\eps})\subset Y_{Q(|\kappa|)+1}$, and continuity of
  $\PRconj|_{Z_{\eps}}$ in 
  $x_0$ follows from the previous lemma.

  To prove that ${\PRconj}|_{Z_{\eps}}$ is a homeomorphism onto its image, it
   remains to show that 
   $\PRconj(x_n)\to\infty$ as $x_n=(\s^n,t_n)\to \infty$ in $Z_{\eps}$. 
   If $s^n_1\to\infty$, then $\im\PRconj(x)\to\infty$ since
   parameter rays cannot
   intersect the lines $\{r+(2k+1)\pi i:r\in\R, k\in\Z\}$, which
   consist entirely of attracting and parabolic parameters. 
   On the other hand, if $(\s,t)\in Z_{\eps}$, then
    \[ T(\F(\s,t)) \geq F(t) - F(t-\eps) \geq
                        (1-\exp(-\eps)) F(t). \] 
   In particular, if $t$ is large enough, then $T(\F(\s,t))\geq t$. 
   By Corollary \ref{cor:parameterraybound}, this implies that 
    $|\PRconj(\s,t)| \geq t/5$
   when $(\s,t)\in Z_{\eps}$ with sufficiently large $t$.  
 \end{proof}

  In particular, we obtain the following analog of Lemma
   \ref{lem:limitset}.

 \begin{cor}[Limit Set of Parameter Rays] \label{cor:limitset}
  Let $\s^0\in\Sequ_0$ and denote the limit
  set of $\PRS_{\s^0}$ by $L$. If there exist some $\s\in\Sequ_0$ 
  and $t>0$ with
  $\PRS_{\s}(t)\in L$, then $\PRS_{\s}\bigl((0,t]\bigr)\subset L$.
 \end{cor}
 \begin{proof} Analogous to Lemma \ref{lem:limitset}. \end{proof}

 \begin{thm}[Cantor Bouquets in Parameter Space]
  There exists a sequence of closed subsets
   $\Xb_k\subset \Xb$ with the following properties:
  \begin{enumerate}
   \item $\Xb_k\subset \Xb_{k+1}$,  
   \item every connected component of $\Xb_k$ is of the form
          $\{\s\}\times [t,\infty)$ for some $t\geq \ts$, 
       \label{item:brush}
   \item The function ${\PRconj}$ is defined on
     $\Xb_k$ and is a homeomorphism onto its image
     $\J_k := \PRconj(\Xb_k)$.  \label{item:homeo}
   \item $\I\subset \bigcup_k \J_k$.
     \label{item:absorbing}
  \end{enumerate}
 \end{thm}
 \begin{proof}
  (Compare \cite[Proof of Theorem 5.4]{markuslassedierk}.)
  Let $\Xb'_n$ denote the set of all $x\in\Xb$ for which
   $\kappa := \PRconj(x)$ is defined and $x\in \F^{-n}(Y_{Q(|\kappa|)+1})$. 
   Then by Lemma \ref{lem:localcontinuity}, 
   $\PRconj : \Xb'_n\to \C$ is defined, injective and a local homeomorphism. 
   Furthermore, for any sequence $(x_n)$ in $\Xb'_n$, clearly
   $\PRconj(x_n)\to\infty$ if and only if $x_n\to\infty$. 

  It follows from the continuity of $\DRconj^{\kappa}(x)$ in
   $\kappa$ and $x$ that
   $\Xb'_n$ is closed. Thus $\PRconj|_{\Xb_n}$ is 
   a homeomorphism onto its image; let 
   $\Xb_n$ be the union of all unbounded connected components
   of $\Xb'_n$. 
    Then $\Xb_n$ satisfies
   (\ref{item:brush}) and (\ref{item:homeo}). 

 It remains
   to establish (\ref{item:absorbing}). Let $\s\in\Sequ_0$.
   We need to
   show that, for every escaping parameter of the form
   $\kappa_0=\PRconj(\s,t_0)$, there is some $n$ such that
   $A := \s\times [t_0,\infty) \subset \Xb_n'$.  

  By \cite[Theorem 3.2]{markusdierk}, there exists
   $t_1 > \ts$ with $\{\s\}\times [t_1,\infty)\subset \Xb_0$.
   Let us set $K := \max_{t\in[t_0,t_1]}|\PRconj(\s,t)|$;
   we can then $n\geq 0$ such that
   $\F^n(A)\subset Y_{K+2}$. Then 
   $(\{\s\}\times [t_0,t_1])\subset \Xb_n$. We thus have
   $A\subset \Xb_n'$, as required.
   \end{proof}

 \begin{cor}[{\cite[Theorem 5.4]{markuslassedierk}}]
  Every path-connected component of $\I$ consists either of a single
   parameter ray
   or of a single parameter ray at a fast address together with an
   escaping landing point.
 \end{cor}
 \begin{proof}
  By Corollary \ref{cor:limitset}, no parameter ray can land at
  a point which is on another parameter ray. Furthermore, by
  \cite[Theorem A.3]{expcombinatorics}, no other parameter ray can
  accumulate
  at the landing point of a fast parameter ray. The 
  claim now follows by applying Proposition 
  \ref{prop:pathcomponentsarecurves},
  where $I:=\I$ and
  $I_n := \J_n\cap \I$ with $\J_n$ from the previous theorem. 
 \end{proof}

\section{Further Questions} \label{sec:openquestions}
 One of the most intriguing question arising from our results is whether
 Theorem \ref{thm:boettcheranalog} generalizes to other classes of entire
 functions, such as the class $\B$ of functions whose set of
 singular values is bounded. Little can be said about the escaping
 sets of such functions without further restrictions. Although for many
 $f\in\B$, the escaping set is arranged in dynamic rays, there also
 exist such functions whose Julia sets contain no curves to $\infty$
 at all \cite{fatoueremenko}. In view of these facts it is perhaps  
 surprising that an analog of Theorem \ref{thm:boettcheranalog}
 \emph{does} hold in full generality for class $\B$
 \cite{boettcher}: if $f,g\in\B$ are
 quasiconformally equivalent in the sense of \cite{alexmisha}, then
 for sufficiently large $R>0$
 there is a conjugacy between $f$ and $g$ defined on the set
 $\{z:|f^n(z)|\geq R\text{ for all $n\geq 0$}\}$. Furthermore,
 this conjugacy extends to a quasiconformal map of the plane. 

In our arguments of escaping set rigidity for Misiurewicz parameters,
 we used the fact that there is an asymptotic value in the Julia set
 which interferes with the topology of the escaping set. 
 The same fact is used in proving 
 the existence of nonlanding rays for exponential Misiurewicz parameters
 (see \cite{nonlanding}). 
 Schleicher \cite{dierkcosine} has shown
 that
 \emph{all} dynamic rays of Misiurewicz members of the
 cosine family $z\mapsto a\exp(z)+b\exp(-z)$ 
 land. 
 It is therefore an
 interesting question whether rigidity of escaping dynamics remains
 for functions without asymptotic values. 

\begin{question}[Escaping Dynamics in the Cosine Family]
 Are there two distinct Misiurewicz parameters in the cosine family
  which are topologically conjugate on their sets of escaping points?
\end{question}

Let us now return to exponential dynamics. As
 mentioned previously, little is known of the accumulation behavior of
 dynamic rays in general. In the case of quadratic polynomials, it is
 still unknown whether a dynamic ray can accumulate on the entire
 Julia set (although it is known \cite{kiwi} that this could happen only
 for Siegel or Cremer parameters). In the exponential family, this
 possibility becomes even more disconcerting:

\begin{question}[Rays Accumulating on the Plane]
 Can the accumulation set of a dynamic ray be the entire complex plane?
\end{question}

We can ask a stronger question:

\begin{question}[Accumulating Rays]
 If $\s^n$ is a sequence of addresses with $|\s_1^n|\to\infty$, is it true
  that $z_n\to\infty$
  whenever $z_n\in g_{\s^n}$ for all $n$?

 More generally, is this true whenever $(\s^n)$ converges to an address
  which is not exponentially bounded?
\end{question}

Let us now depart from questions concerning single rays and
 consider the escaping sets of exponential maps in their entirety. We
 have already formulated the conjecture that two exponential maps
 whose singular value lies in the Julia set are never conjugate on
 their escaping sets by an order-preserving conjugacy. We can ask
 whether the map is already determined by the topology of this
 set. We say that a homeomorphism between
 $I(E_{\kappa_1})$ and $I(E_{\kappa_2})$ is \emph{natural} if it
 preserves the addresses of dynamic rays. 

\begin{question}[Natural Homeomorphisms]
 If $I(E_{\kappa_1})$ and $I(E_{\kappa_2})$ are naturally
 homeomorphic, are $E_{\kappa_1}$ and $E_{\kappa_2}$ conjugate on
 their sets of escaping points?
\end{question} 

The answer to this question can be seen to be ``yes'' when
 $\kappa_1$ and $\kappa_2$ are Misiurewicz-parameters, using the
 construction of nonlanding dynamic rays \cite{nonlanding}. 
 With some more care one can also do this when $\kappa_1$ and
 $\kappa_2$ are escaping parameters lying on different parameter rays.
 The first interesting case in which to investigate thus seems to be that 
 of two parameters on the same parameter ray; for example
 $\kappa_1,\kappa_2\in (-1,\infty)$.

We have described the escaping dynamics completely only in the case of
 attracting and parabolic dynamics. As we
 have seen, the situation becomes much more complicated when the
 singular value moves into the Julia set. Nevertheless, Misiurewicz
 parameters are uniquely determined by their combinatorics. One would
 thus hope that their topological dynamics can also be completely
 understood in terms of their combinatorics, which again might be a
 starting point to understand also more complicated types of
 exponential dynamics.

\begin{question}[Topological Dynamics of Misiurewicz Maps]
 Let $\kappa$ be a Misiurewicz parameter. Can one construct a model
 for the  topological dynamics
 of $\Ek|_{I(\Ek)}$ in terms of $\extaddr(\kappa)$?
\end{question}

Our deliberations in Section \ref{sec:parameterspace} quite naturally
lead to the question of further continuity properties of the map ${\PRconj}$.

\begin{question}[Pinched Cantor Bouquet]
 Is the map ${\PRconj}:X\to \{\kappa:\kappa\in I(\Ek)\}$ a homeomorphism?
  Does ${\PRconj}$ extend to a continuous (and surjective) map from $\Xb$ to
  the exponential bifurcation locus?
\end{question}

A positive answer to the second question would imply density
 of hyperbolicity.

Finally, we have seen that the notion of renormalization fails
 topologically even for attracting parameters. However, there seem 
 to be similarity features in the parameter
 space of exponential maps.
 So it is natural to ask whether some --- different ---
 notion of renormalization might exist in the exponential
 family. Let us formulate this (vaguely) in a special case. Let $W$ be
 any hyperbolic component in exponential parameter space, and let $W'$
 be the unique hyperbolic component of period $1$. Then the multiplier
 maps $\mu:W\to \Ds$ and $\mu':W'\to \Ds$ are universal covering maps
 \cite{alexmisha}, so we can continue some branch of $\mu'^{-1}\circ
 \mu$ to obtain a biholomorphic map 
  $\mathcal{R} : W\to W'$.
 By the results of \cite{habil,boundary}, this map extends
 to a homeomorphism $\mathcal{R}:\cl{W}\to\cl{W'}$.

 \begin{question}[Renormalization]
  Is there some analytic way to construct the parameter
   $\mathcal{R}(\kappa)$ from the dynamics of $\kappa$ in such a way
   that dynamical features such as linearizability etc. are preserved?
 \end{question}

\subsection*{Acknowledgements.} 
  This work was part of my doctoral thesis \cite{thesis}, and I
  would like to thank my advisor, Walter Bergweiler, for all his help
  and encouragement. I would also like to thank Markus F\"orster, Adam
  Epstein, Carsten Petersen, 
  Phil Rippon, Dierk Schleicher and Gwyneth Stallard, 
  as well as the audience
  members at the seminars in which these results were presented for many
  interesting discussions. Furthermore, I am indepted to the
  Institute of Mathematical Sciences at the SUNY Stony Brook and to the
  Mathematics Instute at the University of Warwick
  for continued support and hospitality.
    
\bibliographystyle{hamsplain}
\bibliography{C:/Latex/Biblio/biblio}

\end{document}